\documentclass [pdflatex,sn-mathphys,fleqn]{sn-jnl}

\jyear{2021}

\theoremstyle{thmstyleone}%

\theoremstyle{thmstyletwo}%

\theoremstyle{thmstylethree}%

\usepackage{array}
\usepackage{amsmath}
\usepackage[normalem]{ulem}
\usepackage{etoolbox}
    \makeatletter
    \patchcmd{\ps@headings}
    {\hbox to \hsize{\hfill Springer Nature 2021 \LaTeX\ template\hfill}}
    {\hbox to \hsize{\hfill \hfill}}
    {}
    {}
    \patchcmd{\ps@titlepage}
    {\hbox to \hsize{\hfill Springer Nature 2021 \LaTeX\ template\hfill}}
    {\hbox to \hsize{\hfill \hfill}}
    {}
    {}
    \makeatother

\begin{document}
         
\title[Optimization Framework For Sustainable Post-Disaster Reconstruction]{ A Versatile Optimization Framework For Sustainable Post-Disaster Building Reconstruction}

\author[1]{\fnm{Niloufar} \sur{Izadinia}}
\equalcont{These authors contributed equally to this work.}

\author[2]{\fnm{Elham} \sur{Ramyar}}
\equalcont{These authors contributed equally to this work.}

\author[3]{\fnm{Maytham} \sur{Alzayer}}

\author[3,4]{\fnm{Stephen} \sur{H. Carr}}

\author[2]{\fnm{Gianluca} \sur{Cusatis}}

\author[4,5]{\fnm{Vidushi} \sur{Dwivedi}}

\author[4,6]{\fnm{Daniel} \sur{J. Garcia}}

\author[7]{\fnm{Missaka} \sur{Hettiarachchi}}

\author[1]{\fnm{Thomas} \sur{Massion}}

\author[4]{\fnm{William} \sur{M. Miller}}

\author*[1]{\fnm{Andreas} \sur{W\"achter}}
\email{andreas.waechter@northwestern.edu}

\affil[1]{\orgdiv{Department of Industrial Engineering and Management Sciences}, \orgname{Northwestern University}, \orgaddress{ \city{Evanston}, \state{Illinois}, \country{USA}}}

\affil[2]{\orgdiv{Department of Civil and Environmental Engineering}, \orgname{Northwestern University}, \orgaddress{ \city{Evanston}, \state{Illinois}, \country{USA}}}

\affil[3]{\orgdiv{Department of Materials Science and Engineering}, \orgname{Northwestern University}, \orgaddress{ \city{Evanston}, \state{Illinois}, \country{USA}}}

\affil[4]{\orgdiv{Department of Chemical and Biological Engineering}, \orgname{Northwestern University}, \orgaddress{ \city{Evanston}, \state{Illinois}, \country{USA}}}

\affil[5]{\orgname{Currently at Intel Corporation}, \orgaddress{ \city{Portland}, \state{Oregan}, \country{USA}}}

\affil[6]{\orgname{Currently at Booz Allen Hamilton Inc.}, \orgaddress{ \city{Washington}, \state{DC}, \country{USA}}}

\affil[7]{\orgdiv{Environment and Disaster Management Program},\orgname{World Wildlife Fund}, \orgaddress{ \city{Washington}, \state{DC}, \country{USA}}}

\abstract{This paper proposes an optimization framework for sustainable post-disaster building reconstruction. Based on mathematical optimization, it is intended to provide decision makers with a versatile tool to optimize building designs and to explore the trade-off between costs and environmental impact (represented by embodied energy) of alternative building materials. The mixed-integer nonlinear optimization model includes an analytical building model that considers structural and safety constraints and incorporates regional building codes. Using multi-objective optimization concepts, Pareto-optimal designs are computed that represent the best trade-off designs from which a decision maker can choose when they take additional criteria into consideration. As a case study, we consider the design of a multi-room one-story masonry building in Nepal. We demonstrate how the framework can be employed to address a variety of questions, such as the optimal building design and material selection, the sensitivity of the decision to material prices, and the impact of regional safety regulation thresholds.}


\maketitle

\section{Introduction}\label{sec1}
\label{S:0}

The increase in frequency, intensity, and damage due to extreme weather events is the result of climate change \cite{diffenbaugh2017quantifying}, destruction of coastal ecosystems, deforestation, land degradation, wetland destruction, and stream ecosystem degradation \cite{acreman2013wetlands,thampanya2006coastal}. In addition, the vulnerability to earthquakes and tsunamis may be increased by environmental degradation \cite{van2018analysis,tanaka2007effectiveness}. In 2017, global economic losses from such disasters totaled \$330 billion (USD) \cite{swissre2019}. Furthermore, six of the seven costliest global disasters have occurred within the last ten years, carrying a combined price tag of \$1.7 trillion \cite{swissre2019}. The Intergovernmental Panel on Climate Change (IPCC) warns that, without sweeping reductions in greenhouse gas emissions, global warming will exceed 2°C during the 21st century  \cite{forster2021climate}. This will further increase the frequency of and damage due to hydro-meteorological disasters \cite{field2012managing}. The increasing damage to properties and infrastructure leads to considerable reconstruction investments that result in greater consumption of natural resources \cite{wwfArc2010}.

Therefore, developing an interdisciplinary approach to post-disaster reconstruction for minimizing the costs and adverse environmental effects while keeping resilience and sustainability at acceptable levels is necessary. 
 
\subsection{Post-disaster reconstruction and environmental impact}

Considering the increasing frequency of disasters, minimizing the environmental impact of post-disaster reconstruction is of great importance. Past experiences show that excessive building material extraction (e.g., mining and logging) and processing led to large-scale environmental damage. Hazard-prone and socially inappropriate designs and construction practices, and a shortage of land structurally sound for construction and safe from hazards, were very common among the events studied \cite{vijekumara2015study,jamshed2019building,Shelter2021}. Some of these issues were experienced in the aftermaths of the 2004 Indian Ocean Tsunami and the 2010 Haiti earthquake. Tsumani recovery in Sri-Lanka led to indiscriminate mining of river sand as the primary fine aggregate used in construction, which resulted in riverbank erosion, habitat destruction, riverbed deepening, groundwater depletion, and saltwater intrusion \cite{SriLanka2008}. Disaster debris posed major environmental and social problems in the aftermath of the 2010 Haiti Earthquake \cite{UNDP}. Another side effect of environmental impacts of post-disaster reconstruction is increased greenhouse gas (GHG) emissions that contribute to global climate change and increase vulnerability to future hazards. In addition, incompatibility of existing construction standards, norms, and codes with the needs of post-disaster requirements caused delays in housing reconstruction in many regions in past disasters such as the Indian Ocean Tsunami in 2004.

Practical guidelines for responsible material selection and use for government agencies, nongovernmental organizations (NGOs), and the private sector are rare. The World Wildlife Fund (WWF) and partners developed the Building Material Selection and Use: An Environmental Guide (BMEG) \cite{nepalhariyo} to fill this gap. The BMEG was initially developed for use after the 2015 earthquake in Nepal by WWF Nepal and the Hario Ban Program. Over 700,000 buildings were destroyed or damaged in the earthquake \cite{ministrySci2015} and the demand for common building materials was expected to increase to unprecedented levels. When designing the first version of the BMEG for Nepal, the complex nature of environmental impacts of 26 building materials commonly used in the area was directly addressed. A material’s environmental impacts vary along its own life cycle and the stages of the project. These impacts cannot be minimized simultaneously. However, most of them can be addressed through proper material selection, sourcing, use, and disposal. Finding suitable “environmental intervention points” in the construction project cycle is a key aspect of environmentally responsible disaster reconstruction, also known as the Green Recovery and Reconstruction approach \cite{wwfArc2010}.

Communities' and regions' vulnerability may increase through inappropriate settlement and infrastructure planning and excessive natural resources exploitation. This demonstrates the necessity for an integrated framework for infrastructure planning. 

\subsection{Lack of economic and environmental evaluation}
Few studies have examined the economic and environmental evaluation of disaster reconstruction. Atmaca and Atmaca \cite{atmaca2016comparative}
evaluated the life cycle energy use and costs of temporary housing in Turkey under different scenarios of material costs for container and prefabricated housing. Hosseini et al.\ \cite{hosseini2016multi} used the Integrated Value
Model for Sustainable Assessment (MIVES) that includes the value function concept to evaluate the sustainability value of temporary housing unit (THU) technologies that had been suggested for the 2003 Bam earthquake recovery program. 
One reason that the set of studies is so limited is that post-disaster housing reconstruction environmental impacts are complicated and are most effectively evaluated by research groups including multiple disciplines. The small number of studies demonstrates the need for new optimization tools for minimizing the costs and adverse environmental effects, while keeping projects sustainable and resilient. Such tools would be useful for decision-making after a disaster happens and/or as a preparation guide for potential disaster reconstruction under a set of predicted circumstances.

The remainder of the this paper is organized as follows:  Section~\ref{S:1} reviews previous approaches that used mathematical optimization for sustainable reconstruction, and outlines the envisioned role of the proposed optimization framework. Section~\ref{Sec:BuildingModel} details the mathematical building model used for a one-story masonry building and Section~\ref{Sec:MIPModel} completes the formal optimization problem and the complete mathematical optimization model is presented in the Appendix. 
Section~\ref{sec:results} introduces a case study and describes application of the optimization framework to different use scenarios for the case study. The final discussion in Section~\ref{S:7} includes comparisons with methods used by other investigators.

\section{Background and conceptual approach}
\label{S:1}

Mathematical optimization is widely utilized in all areas of science and engineering. However, only a few case studies have explored its potential to simultaneously minimize the cost and environmental impacts for post-disaster building reconstruction \cite{garcia2015multiobjective, chae2010optimization, venema2005forest}. In this paper, we demonstrate how mathematical optimization can be integrated with a building model that incorporates structural engineering concepts and building codes (Section \ref{Sec:BuildingModel}), and applied to provide rapid and environmentally responsible practical suggestions for post-disaster housing using the example of masonry structures.

\subsection{Mathematical optimization in post-disaster and sustainable reconstruction}

Optimization approaches for evaluating post-disaster and sustainable reconstruction projects are not new. In particular, there are several studies on optimization of multiple objectives.
Dragovic et al.\ \cite{dragovic2017minimiziranje} evaluated the trade-offs between cost and project duration for flood control projects in Serbia. They used the Critical Path Method project scheduling technique and linear programming to solve the optimization problem of minimizing the direct cost of construction subject to constraints such as given deadline, precedence constraints, and upper- and lower-bounded duration time of activities. Their approach predicted significant cost savings for four flood control project case studies. A more recent work by Ghannad et al.\ \cite{ghannad2020multiobjective} proposed a multi-objective optimization model for recovery project prioritization in post-disaster reconstruction to weigh socioeconomic benefits against interruption costs and reconstruction time. For socioeconomic analysis, they used multiple-criteria decision methods to account for social vulnerability to hazards including economics, occupation, ethnicity, etc. They minimized deviation from socioeconomic benefits, along with the reconstruction time and interruption cost, and accordingly found the optimal priority sequence for damaged facilities reconstruction. 

Several other studies focused on debris management and/or transportation networks in a multi-objective or single-objective optimization setting.
Majumder et al.\ \cite{majumder2020review} developed a mathematical optimization model for minimizing the debris transportation cost in a post-earthquake situation and for selecting debris dumping sites along with recycling debris. 
Onan et al.\ \cite{onan2015evolutionary} proposed a multi-objective optimization model for allocating temporary storage facilities for recyclable debris and optimizing debris collection and transportation after a disaster with minimum cost and minimum risk from hazardous waste exposure. El-Anwar et al.\ \cite{el2016efficient} proposed a mixed-integer linear programming model for accelerating transportation network recovery after a disaster. The model is capable of prioritizing different road recovery projects by traffic analysis to minimize the traffic disruption and reconstruction costs. 

Finally, closely related to the approach taken in our paper, Castro-Lacouture et al.\ \cite{castro2009optimization} presented a mixed-integer model for LEED-certified buildings in Colombia. Their model maximizes the number of Leadership in Energy and Environmental Design (LEED) credits subject to design and budget constraints. The authors show that small increases in the budget can lead to a significant increase in credits and present a modified model that determines the minimum budget required to reach a desired number of credits. They also report that the availability of certain building materials is crucial for the achievement of credits, since sustainable materials can be expensive or scarce in a particular region.

Although mathematical optimization has been utilized in post-disaster management, one important component not addressed in the literature to our knowledge is quickly screening multiple possibilities that meet the criteria to develop viable options for post-disaster housing that can be presented to stakeholders. While most of the above studies focused on optimizing overall project sustainability, decision support for specific steps in building construction (e.g.,\ material selection) remains rare.

In this paper, we use bi-objective mathematical programming optimization to generate the best possible combinations of material selection and construction decisions for post-disaster housing to minimize the required budget and embodied energy. Embodied energy is a measure of a resource requirement that a material represents when it is made into a useful form, such as a building block.  It is an intrinsic number calculated by summing all the energy inputs, from obtaining raw material and continuing through refining and finishing steps. Transportation costs are also included in this summation. We use embodied energy as one measure of environmental impact that is correlated with greenhouse gas emissions.

\subsection{Envisioned role of the optimization tool in the post-disaster reconstruction process}

\begin{figure}[ht]
\begin{center}
\includegraphics[width=1\textwidth]{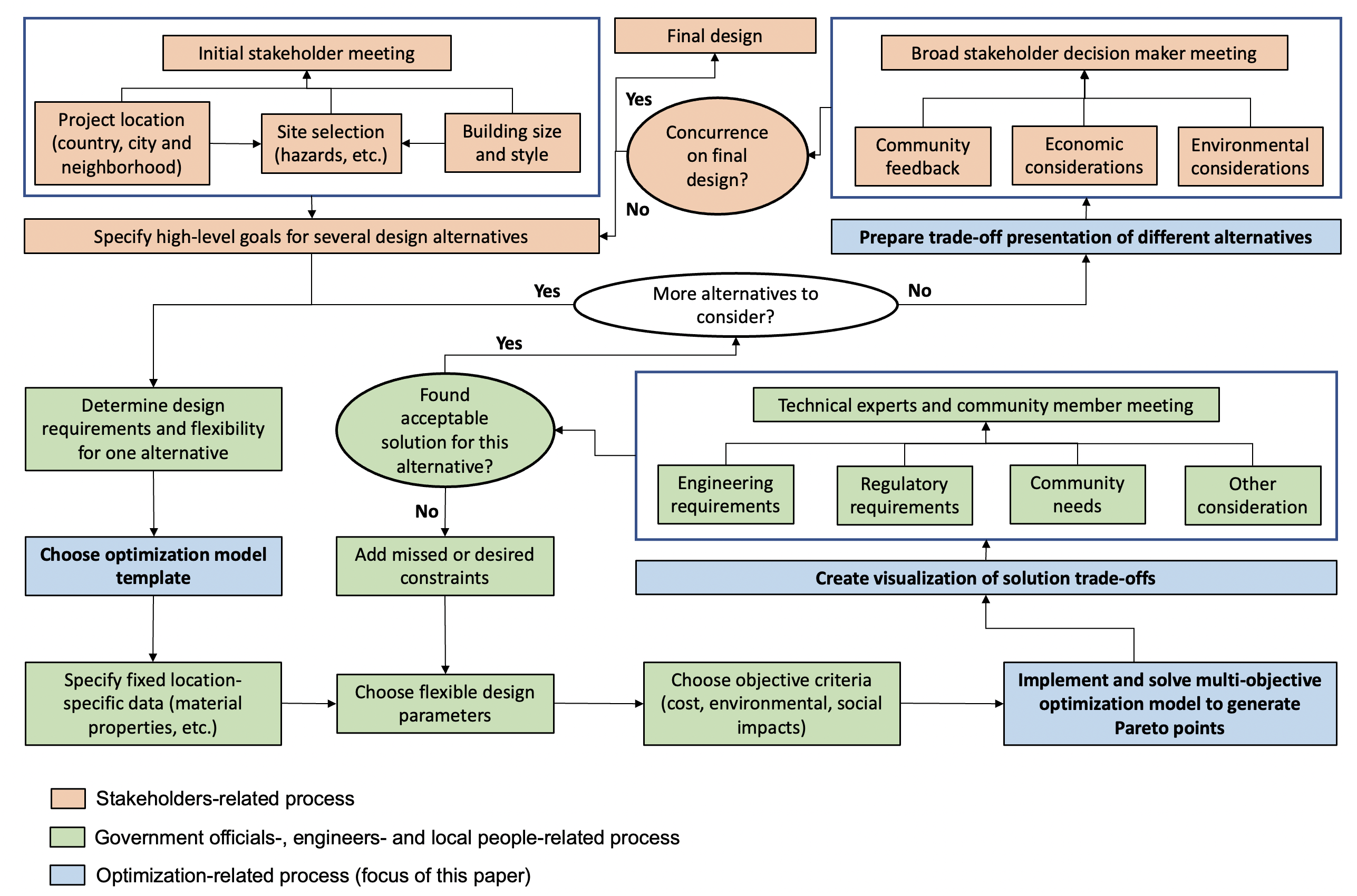}
\end{center}
\caption{Master planner for development of alternatives and selection of final building design process flowchart}
\label{fig:flowchart}
\end{figure}

We aim to use optimization that is integrated with a building model to rapidly develop alternative design solutions that meet all structural requirements and regulations. Our goal is to develop an interactive and iterative system that engages stakeholders at the beginning of the process and helps them visualize the best trade-off between cost and environmental impacts. The envisioned process (Figure \ref{fig:flowchart}) begins with an initial stakeholder (e.g., government officials, NGOs, and funders) meeting to determine high-level goals including the project location, site selection, and building size and style for several design alternatives. Then, the experts (e.g., engineers, architects, procurement manager, and project manager) determine design requirements and flexibility for each alternative. They provide this information to the optimization team to determine the appropriate mathematical optimization model. The initial model and material properties in the database are adapted in terms of site- and design-specific constraints, and input data and flexibility, according to the stakeholders’ and experts’ desire. In addition, the objective criteria are chosen for the model. In this paper, the criteria relate to cost and embodied energy, as a measure of environmental impact. However, other options such as water consumption or social impact are possible. The adapted model is solved for different optimal alternatives along with visualization of the trade-off between objectives. Then, a meeting with technical experts and community members is held to discuss the results of the model, after which changes to the model are made if necessary. The modified model is solved and the iterative process continues until a solution acceptable to experts and community members is obtained. These results are presented to the broader stakeholders for their feedback. The outer iterative process continues until the design is accepted by community members, experts, and broader stakeholders.

To demonstrate the range of questions that the proposed optimization methodology can address, this paper incorporates architectural considerations to obtain an optimal solution minimizing the cost and environmental impact in a multi-room one-story masonry building. The process in Figure~\ref{fig:flowchart} is a generic master planner which may be applied to many situations. 

Our approach is designed to optimize permanent post-disaster reconstruction, rather than providing temporary shelters. This also differs from long-term planning due to the shorter time-frame and need to build many structures at the same time.

\section{Building model} 
\label{Sec:BuildingModel}

The purpose of the proposed optimization framework is to quickly identify optimal building designs while ensuring that the resulting structures are safe, meet local building codes, and satisfy residents’ demands for space and comfort. Therefore, to account for structural engineering concepts, mechanics and materials, and design codes, a building model must be integrated into the optimization framework. The building model is critical to ensure that the optimized building design and material selections are resilient and able to withstand operation loads, load combinations, and future disasters. 
In this section, we explain the constructional constraints in the mathematical model. We postpone introducing the objective functions to Section \ref{Sec:MIPModel}.

\begin{figure*}[ht]
	\centering
	\includegraphics[width=0.75\textwidth]{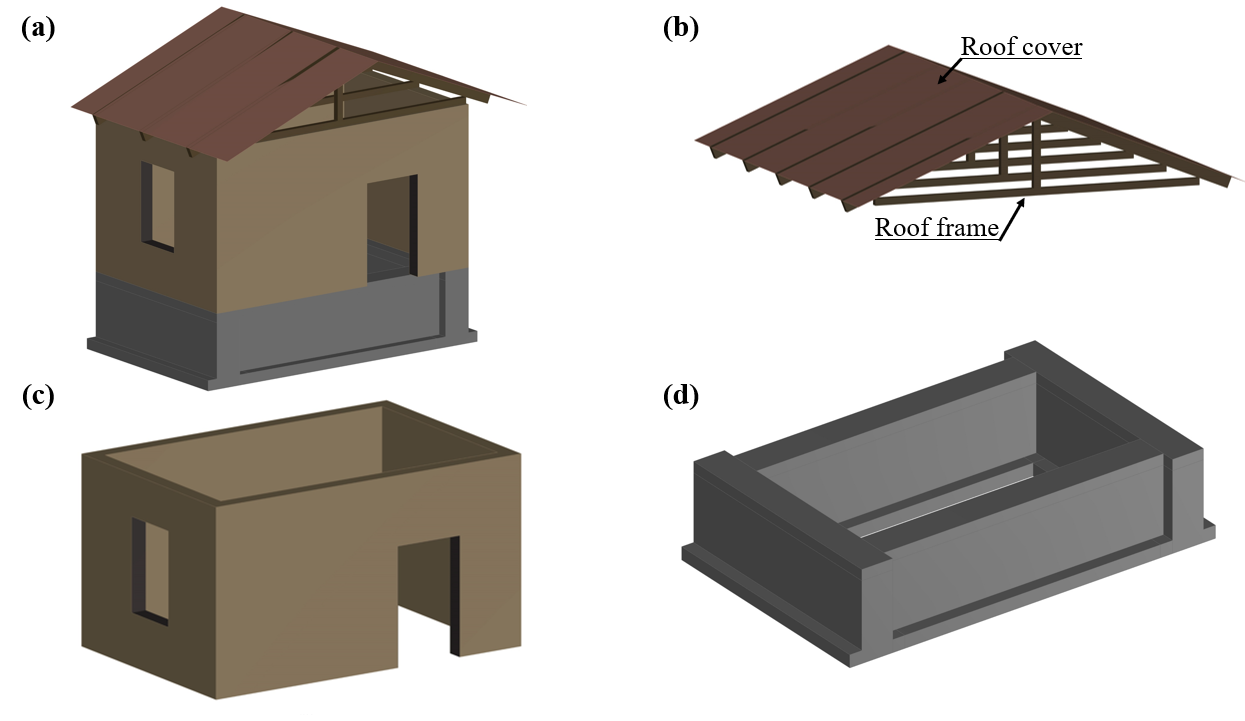}
	\caption{Typical one-story masonry building (a) perspective view; (b) roof; (c) superstructure; (d) foundation}
   \label{fig:one-story-masonry-house}
\end{figure*}

\begin{figure*}[ht]
	\centering
	\includegraphics[width=0.75\textwidth]{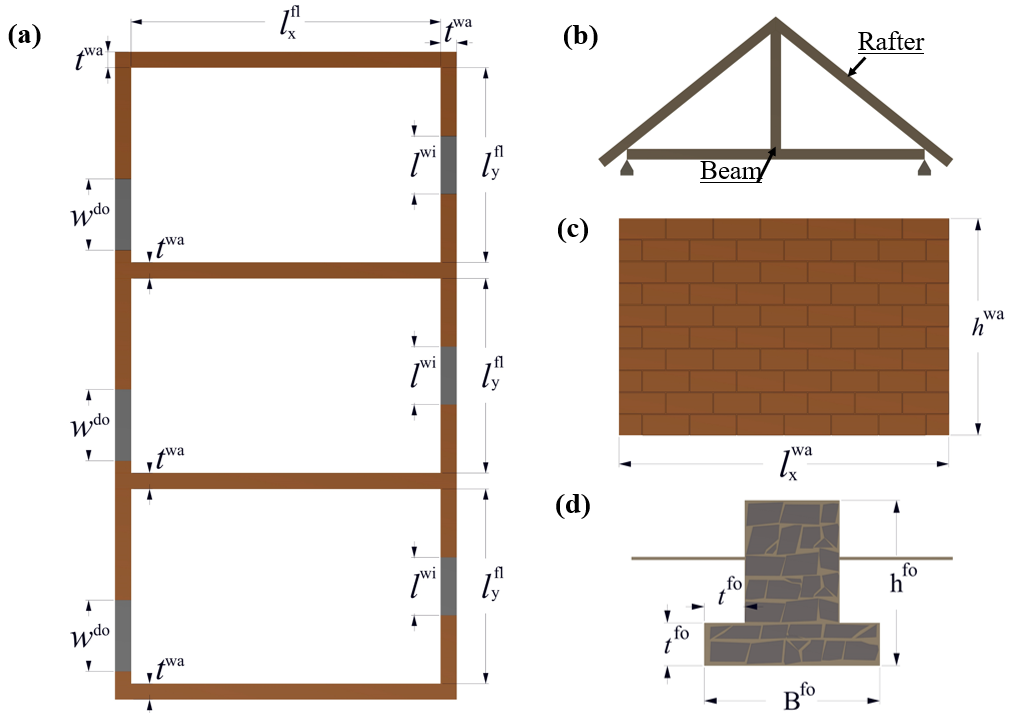}
	\caption{Masonry building (a) typical floor plan view of a three-room school (or home) used for this analysis; (b) typical elements of one slice of the roof frame; (c) typical brick wall; (d) typical stone foundation}
   \label{fig:geometry-masonry-house}
\end{figure*}

To this end, the mathematical optimization model integrates an analytical building model to determine material shapes and amounts for the main structural elements (i.e., roof, superstructure, and foundation shown in Figures \ref{fig:one-story-masonry-house} and \ref{fig:geometry-masonry-house}). For demonstration purposes, we consider multi-room masonry buildings. The building model includes the design criteria and restrictions of the Nepal National Building Codes \cite{dudbc} such as \cite{Unreinforced, Seismic, BearingLoad, Occupancy, Earthquake, Infill}. In the absence of information from these sources, construction data from India’s sample sites \cite{ MortarPaste, FiredBrick1, FiredBrick2, CompressedSoil1, CompressedSoil2, UnitofMasonary1, UnitofMasonary2} are used.
We stress that the optimization framework proposed in this paper is not limited to a particular building model and could be applied to building types other than multi-room masonry buildings and different local building codes, and could also consider other optimization criteria.

The model mimics the construction criteria for a one-story building with multiple adjacent rooms that could be used as a clinic or school, etc. (Figure \ref{fig:geometry-masonry-house}). Therefore, the layout of the building is considered to be a row of rooms that have common walls. In the remainder of this section, the formulated building model is presented in detail. Table~\ref{table:sets} specifies the sets of available building materials. Tables~\ref{table:parameter} and \ref{table:variables} define the notation of the parameters and variables pertaining to the roof, superstructure, and foundation of the model. 

These quantities are explained in detail below, along with different sets of structural feasibility constraints. Figure~\ref{fig:geometry-masonry-house} shows the roof frame elements, and dimension parameters of the floor plan of a typical three-room school (or home), brick wall, and stone foundation of a typical masonry building. Table~\ref{table:parameter} listed the parameters, i.e., the quantities with values that are fixed for the optimization. The parameters generally denoted by roman letters (such as ${\rm C_m}$), while the variables, i.e., the quantities that are adjusted during the optimization, are generally denoted by italic letters (such as $v^{\text{wa}}$). Table~\ref{table:variables} shows the variables. 
Table~\ref{table:parameter} also gives the parameter values used in the case study. For each optimization variable, Table~\ref{table:variables} shows its feasible range or gives the equation by which the variable value is computed. 
The full optimization model is given in Appendix~\ref{secmodel}.

\begin{table}[ht]
 \begin{center}
 \begin{minipage}{250pt}
\caption{Sets of available building materials.}\label{table:sets}
\begin{tabular}{@{}lll@{}}
  \toprule
	Set & Component & Options \\
     \midrule
     $c \in \mathcal{C}$ & roof covering  & \{bamboo, plywood\}  \\
     $f \in \mathcal{F}$ & foundation  & \{stone, brick, concrete block, soil block\}  \\ 
    $r \in \mathcal{R}$ & roof  & \{bamboo, wood\} \\
    $w \in \mathcal{W}$ & wall  & \{stone, brick, concrete block, soil block\} \\
    $m \in \mathcal{M}$ & all materials & $\mathcal{C}\cup\mathcal{F}\cup\mathcal{R}\cup\mathcal{W}$ \\
\bottomrule
\end{tabular}
\end{minipage}
\end{center}
\end{table}

\begin{table*}
  \begin{scrisize}
 \caption{Parameters in the building and optimization model, including their values.} \label{table:parameter}
 \begin{center}         
 \begin{tabular}{@{}p{0.08\textwidth} p{0.45\textwidth} p{0.08\textwidth}p{0.22\textwidth} @{}}
    \toprule
          \multicolumn{1}{c}{\textbf{Parameter}} & \multicolumn{1}{c}{\textbf{Explanation}} &
          \multicolumn{1}{c}{\textbf{Units}} &
          \multicolumn{1}{c}{\textbf{Value}} \\
        \cmidrule{1-4}
    \multicolumn{4}{l}{\textbf{General parameters:}}  \\
    ${\rm B^{\text{avail}}}$ & available budget  &\$& varies \\
        ${\rm C}_m$  & unit cost of material $m \in  \mathcal{M}$ &\$/m$^3$& varies  \\
    ${\rm E}_m$  & unit embodied energy of material  $m \in  \mathcal{M}$  &MJ/kg&  varies   \\
 ${\rm n^{\text{rm}}}$  & number of rooms &-& 3  \\
     $\rho_m$  & unit density of material   $m \in  \mathcal{M}$ &kg/m$^3$& varies \\
     {\rm g}&acceleration of gravity&N/kg&9.8\\
    \cmidrule{1-4}
      \multicolumn{4}{l}{\textbf{Superstructure (floor, wall, and openings) parameters:}}   \\
    ${\rm \underline{A}^\text{fl}}$ & minimum floor space area per room &m$^2$& 10  \\
    ${\rm \sigma}^\text{allw,c}_{w}$ & allowable compressive strength of wall material $w\in\mathcal W$ &MPa& varies \\
    ${\rm C^\text{re}}$  & unit cost of rebar material  &\$/m& 17.3 \\
    ${\rm C_\text{f}}$ & wind load coefficient &-& 1 \\
    $\{\rm C,\rm Z,\rm I,\rm K\}$ & seismic load components &-& \{0.08, 1.1, 1, 4\} \\
    ${\rm C_\text{d}}$ & seismic load coefficient &-& ${\rm C_\text{d}} = \rm C \cdot \rm Z \cdot \rm I \cdot \rm K$ \\ 
    ${\rm d^\text{re}}$  & rebar diameter &m& 0.012 \\
    ${\rm E^\text{re}}$  & unit embodied energy of rebar material  &MJ/kg&  37.95 \\
    $\{{\rm \underline{h}^\text{wa}},{\rm \bar{h}^\text{wa}}\}$ & minimum and maximum wall height &m& \{2.7, 3.8\} \\
    ${\rm h^\text{do}}$ & door height &m& 2 \\
    ${\rm \underline{l}^\text{wi}}$ & minimum window length &m& 0.7 \\
    ${\rm \underline{n}^\text{re}}$  & minimum number of rebar slices &-& 2 \\
    ${\rm P_\text{L}}$ & live load per square meter &$\text{N}/\text{m}^2$& 2,000 \\
    ${\rm P^\text{design}}$ & weight applied by wind on the structure  &   $\text{N}/\text{m}^2$ & 8,120\\
    $\rho^\text{re}$  & linear density of rebar material & kg/m&  0.89  \\
    ${\rm \tau^\text{allw}}$ & allowable shear &MPa& 0.5 \\
    ${\rm \underline{s}^\text{re}}$  & minimum inter-rebar spacing &m& 0.05 \\
    ${\rm \underline{t}^\text{wa}_{w}}$ & minimum wall thickness for $w \in \mathcal{W}$&m& varies  \\
    ${\rm \bar{t}^\text{wa}}$ & maximum wall thickness &m& 1.1  \\
    ${\rm \sigma^\text{allw}_{t}}$ & allowable tension &MPa& 0.12 \\
    ${\rm \underline{w}^\text{do}}$ & minimum door width &m& 1.1\\
            \cmidrule{1-4}
          \multicolumn{4}{l}{\textbf{Roof and roof covering parameters:}}  \\
    ${\rm A^\text{be}}$ & cross-sectional area of beam &$\text{m}^2$& 0.0253  \\
    ${\rm A^\text{ra}}$ & cross-sectional area of rafter &$\text{m}^2$& 0.0065  \\
    $\{{\rm \underline{l}^\text{fl}},{\rm \bar{l}^\text{fl}}\}$ & minimum and maximum floor space length &m& \{2, 4.5\} \\
    $\{{\rm \underline{n}^\text{slc}},{\rm \bar{n}^\text{slc}}\}$ & minimum and maximum number of slices &-& \{2, 20\} \\
    ${\rm R^{\text{be}}}$ & ratio of rafter length to beam length &-& 0.4027 \\
    ${\rm R^{\text{co}}}$ & ratio of roof cover material volume to roof material volume &-& 0.0628 \\
    ${\rm \bar{s}^\text{be}}$ & maximum spacing between beams &m& 0.5  \\
    ${\rm w^\text{be}}$ & cross-sectional beam width &m& 0.11 \\
    \cmidrule{1-4}
    \multicolumn{4}{l}{\textbf{Foundation parameters:}} \\
    ${\rm h^\text{fo}}$ & height of foundation &m& $1.1 $  \\
    ${\rm \underline{t}^\text{fo}_{f}}$ & minimum foundation thickness for material $f \in \mathcal{F}$  &m& varies   \\
    ${\rm B^\text{fo}}$ & width of foundation &m& $0.8$ \\
        \bottomrule
     \end{tabular} 
\end{center}
  \end{scrisize}
 \end{table*}

 \begin{table}
 \begin{center} 
                 \caption{The variables in the building and optimization model, together with their physical definition or allowable ranges.  Some of the equations and ranges depend on the building material chosen.  In this case, $w^*\in\mathcal W$ denotes the chosen wall material, and $a^*$, $f^*$, and $r^*$ are similarly defined.}\label{table:variables}
\resizebox{\columnwidth}{!}{
\begin{tabular}{@{} p{0.08\textwidth} p{0.45\textwidth} p{0.03\textwidth} p{0.52\textwidth} @{}}
    \toprule
    \multicolumn{1}{c}{\textbf{Variable}} & \multicolumn{1}{c}{\textbf{Explanation}} &
    \multicolumn{1}{c}{\textbf{Units}} &
    \multicolumn{1}{c}{\textbf{Definitions and restrictions}} \\
       \cmidrule{1-4}
      \multicolumn{4}{l}{\textbf{Superstructure (floor, wall, and openings) variables:}}   \\
$\{A^\text{wa}_\text{x},A^\text{wa}_\text{y}\}$ & area of wall cross section  in $x$- and $y$-direction &m$^2$& $\{t^\text{wa}\cdot l_\text{x}^\text{wa},t^\text{wa}\cdot l_\text{y}^\text{wa}\}$ \\[2ex]
    $F_1$ & seismic load external force &N& $F_1 = {\rm g} \, {\rm C_\text{d}} \, v^\text{wa}\, \rho_{w^*} $ \\[2ex]
    $F_\text{e}$ & seismic force in $x$- or $y$-direction imposed on each wall &N & $F_\text{e} = 0.5\, F_1$ \\[2ex]
    $\{F_{w,\text{x}},F_{w,\text{y}}\}$ & wind force in $x$- and $y$-direction imposed on each wall &N& $0.5\, {\rm C_\text{f}} \cdot \{l_\text{y}^\text{wa},l_\text{x}^\text{wa}\} \cdot h^\text{wa}\cdot {\rm P^\text{design}}$ \\[2ex]
    $ h^\text{wa}$ & wall height  &m&  ${\rm \underline{h}^\text{wa}} \leq$ $  h^\text{wa} $  $\leq {\rm \bar{h}^\text{wa}}$ \\[2ex]
    $ \{l_\text{x}^\text{fl},l_\text{y}^\text{fl}\}$ & floor space length in $x$- and $y$-direction per room &m& ${\rm \underline{l}^\text{fl}}\leq l_\text{x}^\text{fl},\; l_\text{y}^\text{fl} \leq {\rm \bar{l}^\text{fl}}$ \\[2ex]
    $ \{l_\text{x}^\text{wa},l_\text{y}^\text{wa}\}$ & wall length in $x$- and $y$-direction per room &m& $\{  l_\text{x}^\text{fl} + 2\, t^\text{wa}, l_\text{y}^\text{fl}  + 2\, t^\text{wa}\}$ \\[2ex]
    $ l^\text{wi}$ & window length &m& $l^\text{wi} \geq {\rm \underline{l}^\text{wi}}$ \\[2ex]
    $ l^{\text{re}}$ & total length of rebar per room &m& $l^{\text{re}} = n^\text{re} \cdot 2\, ( w^\text{do} + {\rm h^\text{do}} ) +4\, l^\text{wi} $ \\[2ex]
    $\{M_{w,\text{x}},M_{w,\text{y}}\}$ & moment of wind load applied to wall in $x$- and $y$-direction &N.m& $\{F_{w,\text{x}}\, h^\text{wa}, F_{w,\text{y}}\, h^\text{wa}\}$ \\[2ex]
       $M_\text{e,x}$ & moment of seismic load applied to wall in $x$-direction &N.m& $M_\text{e,x} = F_\text{e}\, h^\text{wa}$ \\[2ex]
    $\{S_{\text{x}},S_{\text{y}}\}$ & section modulus of wall in $x$- and $y$-direction &m$^3$& $\{ \frac{1}{6} t^\text{wa} \cdot (l_\text{x}^\text{wa})^2,  \frac{1}{6} t^\text{wa} \cdot (l_\text{y}^\text{wa})^2\}$ \\ [2ex]
    $n^\text{re}$ & number of rebar slices for the doors&-&  $ n^\text{re} \geq {\rm \underline{n}^\text{re}}  ,\; n^\text{re} \in \mathbb{Z}_+$ (positive integer)\\[2ex]
    $P_\text{D}^\text{wa}$ & self-weight of wall per unit length &N/m& $P_\text{D}^\text{wa} = {\rm g}\, t^\text{wa} h^\text{wa}\rho_{w^*}$ \\[2ex]
    $F_\text{D,x}^\text{wa}$ & total dead-load applied to wall in $x$-direction &N& $F_\text{D,x}^\text{wa} = P_\text{D}^\text{wa} l_\text{x}^\text{wa} + 0.5{\rm g}\,q^\text{ro}$  \\[2ex]
    $F_\text{L,x}^\text{wa}$ & total live-load applied to wall in $x$-direction &N& $F_\text{L,x}^\text{wa} = 0.5 {\rm P_\text{L}}\, l_\text{x}^\text{wa}\, t^\text{wa}$ \\[2ex]
    $t^\text{wa}$ & wall thickness &m& ${\rm \underline{t}^\text{wa}_{w^*}}\leq $ $ t^\text{wa} $ $\leq {\rm \bar{t}^\text{wa}}$ \\[2ex]
    $v^\text{wa}$ & total volume of wall material &m$^3$& $v^\text{wa} = t^\text{wa} \, 
     \big(2\, h^\text{wa} \, (l_\text{y}^\text{fl}+l_\text{x}^\text{wa}) - (w^\text{do} \, {\rm h^\text{do}} + (l^{\text{wi}})^2) \big)$ \\[2ex]
  $v^\text{wa}_\text{tot}$ & total volume of wall material for entire building &m$^3$& $v^\text{wa}_\text{tot} = {\rm n^{\text{rm}}} v^\text{wa} - ({\rm n^{\text{rm}}}-1)\, t^\text{wa} h^\text{wa} l_\text{x}^\text{wa} $ \\[2ex]
  $l^\text{re}_\text{tot}$&total length of rebars for entire building &m&$l^\text{re}_\text{tot} = {\rm n^{\text{rm}}} l^{\text{re}}$\\[2ex]
  $ w^\text{do}$ & door width &m& $w^\text{do} \geq {\rm \underline{w}^\text{do}}$ \\[2ex]
    \cmidrule{1-4}
    \multicolumn{4}{l}{\textbf{Roof and roof covering variables:}}\\[2ex]
    $ n^\text{slc}$ & number of roof slices that make up the roof &-& ${\rm \underline{n}^\text{slc}} \leq $ $  n^\text{slc}$ $\leq {\rm \bar{n}^\text{slc}} $ , $ n^\text{slc}$ $\in \mathbb{Z}_+$ \\[2ex]
     $q^\text{ro}$ & total roof weight &kg& $q^\text{ro} = n^\text{slc} \cdot v^\text{slc} (\rho_{r^*} +  {\rm R^{\text{co}}} \rho_{c^*})$ \\[2ex]
    $v^\text{slc}$ & volume of a single slice of roof material &m$^3$& $  v^\text{slc} $ $= {\rm A^\text{be}} l_\text{x}^\text{wa} + 2 {\rm A^\text{ra}} \sqrt{ \big(\frac{{\rm R^{\text{be}}}\cdot l_\text{x}^\text{wa}}{2} \big)^2 + \big(\frac{l_\text{x}^\text{wa}}{2}\big)^2}$ \\[2ex]
    $  v^\text{slc}_\text{tot}$ & total volume of roof material for entire building &m$^3$& $  v^\text{slc}_\text{tot} = {\rm n^{\text{rm}}} \cdot n^\text{slc}$ $\cdot  v^\text{slc}$ \\[2ex]
    $  v^\text{co}_\text{tot}$ & total volume of roof covering material  for entire building &m$^3$& $v^\text{co}_\text{tot}= {\rm R^{\text{co}}}\cdot v^\text{slc}_\text{tot}$ \\[2ex]
    \cmidrule{1-4}
    \multicolumn{4}{l}{\textbf{Foundation variables:}}\\[2ex]
    $e$ & eccentricity &m& $e = 0.5 \cdot ({\rm B^\text{fo}}-2 t^\text{fo} - t^\text{wa})$ \\[2ex]
        $P^\text{fo}_\text{D}$ & self-weight of foundation per unit of length &N/m& $P^\text{fo}_\text{D} = {\rm g}\cdot A^\text{fo} \cdot \rho_{f^*}$ \\[2ex]
    $F^\text{fo}_\text{D,x}$ & total weight laid on $x$-direction of foundation &N& $F^\text{fo}_\text{D,x} = P^\text{fo}_\text{D} l_\text{x}^\text{wa} + F_\text{D,x}^\text{wa}$ \\[2ex]
    $t^\text{fo}$ & foundation thickness &m& $ {\rm \underline{t}^\text{fo}_{f^*}} \leq t^\text{fo} \leq 0.5 {\rm B^\text{fo}}$ \\[2ex]
    $A^\text{fo}$ & cross sectional area of foundation &m$^2$& $A^\text{fo} = {\rm B^\text{fo}} {\rm h^\text{fo}} - 2 t^\text{fo} ({\rm h^\text{fo}}-t^\text{fo})$ \\ [2ex]
    $v^\text{fo}_\text{tot}$ & total volume of foundation material for entire building &m$^3$& $v^\text{fo}_\text{tot} = (2{\rm n^{\text{rm}}}(l_\text{y}^\text{wa}-2t^\text{fo})\; + $  $({\rm n^{\text{rm}}} + 1)l_\text{x}^\text{wa})\cdot A^\text{fo} $ \\ 
    \bottomrule
     \end{tabular} }
\end{center}
 \end{table}

\subsection{Roof constraints}
\label{S:roof_constraints}
The roof of a standard single-story masonry building consists of a collection of triangular slices, each made up of two rafters connecting at the top and a beam at the bottom. These triangular slices are then laid across the walls in the $y$-direction, and then covered with auxiliary materials $\mathcal{C}$ shown in Table \ref{table:sets}. 
The following constraint ensures that the beams fit within the total length $l_\text{y}^\text{wa}$ of the building in the $y$-direction and that the space between the slices does not exceed ${\rm \bar{s}^\text{be}}$, a restriction imposed by the building code.

\begin{equation}\label{1lp}
n^\text{slc} {\rm w^\text{be}} \leq l_\text{y}^\text{wa} \leq  (n^\text{slc}-1) \cdot {\rm \bar{s}^\text{be}} + n^\text{slc} {\rm w^\text{be}} 
\end{equation}

\subsection{Wall constraints}
\label{S:5_1_3}
A key input parameter for the model is the minimum floor space per room of the building (${\rm \underline{A}^\text{fl}}$). This is often determined by those stakeholders who will ultimately use the building.
The minimum floor space is maintained by the following constraint: 
\begin{equation}
    l_\text{x}^\text{fl} \cdot l_\text{y}^\text{fl} \geq {\rm \underline{A}^\text{fl}}.
\end{equation}

Another major limitation from an architectural standpoint is the requirement to withstand various disasters. 
The model considers different applied loads such as self-weight (dead), live loads such as snow and wind, and seismic load. The classification of construction sites and environmental conditions (e.g., bed soil type, earthquake severity, and wind speed)\ are taken into account. Different load types can act on the structure simultaneously, and the model ensures that the actual stresses do not exceed the worst-scenario thresholds.
The limiting load combination inequalities below consider the $x$-direction of the building.  Their analogues for the $y$-direction are also imposed.

According to the basis of mechanics of materials  and allowable stress design criteria \cite{beer1999mechanics}, the shear stress is limited by
\begin{align*}
\frac{3 F_{w,\text{x}}}{2   A^\text{wa}_\text{x}}&\leq {\rm \tau^\text{allw}} \\[2ex]
\frac{3 F_\text{e}}{2 A^\text{wa}_\text{x}}&\leq {\rm \tau^\text{allw}},
\end{align*}
the  tensile stress by
\begin{align*}
 -\frac{F_\text{D,x}^\text{wa}+F_\text{L,x}^\text{wa}}{A^\text{wa}_\text{x}}+\frac{M_{w,\text{x}}}{S_{\text{x}}} &\leq  {\rm \sigma^\text{allw}_{t}} \\
 -\frac{F_\text{D,x}^\text{wa}+F_\text{L,x}^\text{wa}}{A^\text{wa}_\text{x}}+\frac{M_\text{e,x}}{S_{\text{x}}} &\leq {\rm \sigma^\text{allw}_{t}},
\end{align*}
and the cross-sectional compressive stress is limited by
\begin{align*}
\frac{F_\text{D,x}^\text{wa}+F_\text{L,x}^\text{wa}}{A^\text{wa}_\text{x}}+\frac{M_{w,\text{x}}}{S_{\text{x}}} &\leq  {\rm \sigma}^\text{allw,c}_{w^*} \quad\\
\frac{F_\text{D,x}^\text{wa}+F_\text{L,x}^\text{wa}}{A^\text{wa}_\text{x}}+\frac{M_\text{e,x}}{S_{\text{x}}} &\leq {\rm \sigma}^\text{allw,c}_{w^*},
\end{align*}
where $w^*\in\mathcal{W}$ is the selected wall material.

\subsection{Foundation constraints}
\label{S:foundation_constraints}

Similar to the walls, the building code requires that the shear stress experienced by the foundation does not exceed a threshold.  
The formula to compute the stress changes depending on whether the value of the eccentricity $e$, a measure of circularity, is below or above a threshold based on the width of the foundation \cite{coduto2001foundation}.
More specifically, if $e \geq \frac16{{\rm B^\text{fo}}}$, we need to enforce
\begin{align*}
 \left(1 + \frac{t^\text{fo}}{{\rm B^\text{fo}}} \right)\cdot \frac{F^\text{fo}_\text{D,x}+F_\text{L,x}^\text{wa}}{l_\text{x}^\text{wa} ({\rm B^\text{fo}}-2 e)} & \leq {\rm \tau^\text{allw}} \\ 
  \left(1 + \frac{t^\text{fo}}{{\rm B^\text{fo}}}\right) \cdot \frac{P^\text{fo}_\text{D} + P_\text{D}^\text{wa}}{{\rm B^\text{fo}}- 2 e} &\leq {\rm \tau^\text{allw}},
\end{align*}
where the first constraint considers the $x$-direction, and the second the $y$-direction.
On the other hand, when $e < \frac16{\rm B^\text{fo}}$, we need to enforce
\begin{align*}
    \frac{3}{2 \, l_\text{x}^\text{wa}}\cdot(F^\text{fo}_\text{D,x}+F_\text{L,x}^\text{wa}) \cdot \left(\frac{1}{{\rm B^\text{fo}}} + \frac{6 \, t^\text{fo} \, e}{({\rm B^\text{fo}})^3} \right) & \leq {\rm \tau^\text{allw}}\\
    \frac{3}{2} \cdot (P^\text{fo}_\text{D} + P_\text{D}^\text{wa})\cdot \left(\frac{1}{{\rm B^\text{fo}}} + \frac{6 \, t^\text{fo} \, e}{({\rm B^\text{fo}})^3} \right) &\leq {\rm \tau^\text{allw}},
\end{align*}
again for the $x$- and the $y$-direction, respectively.

\subsection{Opening constraints}
\label{S:opening_constraints}

Another set of safety constraints limits the size of the windows and doors of the structure based on the relative size of the wall lengths.
We should ensure that the window length and door width does not exceed half of the larger of the wall length in the $x$- and the $y$-directions (the door and window can be considered into either the $x$- or the $y$-direction wall, however both of them should not be in the same wall): 

\begin{align}\label{eq:openings1:raw}
w^\text{do} &\leq \tfrac{1}{2} \cdot \max\{l_\text{x}^\text{wa},l_\text{y}^\text{wa} \} \\
	l^\text{wi} &\leq \tfrac{1}{2} \cdot \max\{l_\text{x}^\text{wa},l_\text{y}^\text{wa} \} .\label{eq:openings2:raw}
\end{align}

Furthermore,  the window length cannot be larger than half of the wall height:
\begin{align}\label{2lp}
&	l^\text{wi} \leq \frac{h^\text{wa}}{2}.
\end{align}
Note that windows are square and the length and width are the same. 
To ensure the structural integrity of the windows and doors, reinforcement barring, or rebar, is necessary. The following constraint maintains the minimum spacing between rebar frames:
\begin{equation}\label{3lp}
n^\text{re} {\rm d^\text{re}} + (n^\text{re} - 1) {\rm \underline{s}^\text{re}} \leq t^\text{wa}.
\end{equation}

\section{The mixed-integer nonlinear optimization formulation}
\label{Sec:MIPModel}

This section presents a mathematical optimization model that strives to strike a balance between the total embodied energy of the construction and the required budget. 
This model is based on the building model introduced in the previous section as its main constraints.  In particular, the entities listed in Table~\ref{table:variables} are the decision variables that need to be determined by the optimization solver.  Similarly, the definitions and restrictions in Table~\ref{table:variables} and the structural constraints given in Sections~\ref{S:roof_constraints}--\ref{S:opening_constraints} comprise the constraints of the optimization problem.

In addition to the variables listed in Table~\ref{table:variables}, Table~\ref{table:linear} introduces the binary variables $x_m$ that indicate if building material $m\in\mathcal M$ is selected for the respective components.  For example, for each choice $w\in\mathcal W$ of the walls material, the variable $x_w$ is set to 1 if the material is selected and 0 otherwise.  The constraint $\sum_{w\in \mathcal W} x_w=1$ ensures that the optimization model selects exactly one wall material.
Further binary variables are included to express some logical relationships, such as conditional constraints.  These are discussed in more detail in Section \ref{S:logical}.

The overall optimization model is a mixed-integer nonlinear nonconvex optimization problem; note that all constraints, except \eqref{1lp}, \eqref{2lp}, and \eqref{3lp}, are nonlinear and nonconvex. This model can be solved with specialized optimization solvers, such as BARON \cite{sahinidis:baron:21.1.13} and SCIP  \cite{GamrathEtal2020OO}.
The complete optimization model is given in Appendix \ref{secmodel}.

\begin{table}
 \caption{The added decision variables and accordingly modified equations in the mixed-integer nonlinear optimization problem formulation. }\label{table:linear}
 \begin{center}
\begin{tabular}{@{}lll@{}}
  \toprule
    \multicolumn{1}{c}{Variables} & \multicolumn{1}{c}{Explanation} & \multicolumn{1}{c}{Definitions and restrictions} \\
     \cmidrule{1-3}
     \multicolumn{3}{l}{\textbf{Binary variables:}}\\
    $ x_c$ &$=1$ if roof covering material $c\in\mathcal C$ is selected & $\sum_{c \in \mathcal{C}}$ $  x_c $ $ = 1$, $ x_c$ $ \in \{0,1\}$ \\
    $x_f$ & $=1$ if foundation material $f\in\mathcal F$ is selected & $\sum_{f \in \mathcal{F}} x_f  = 1, x_f \in \{0,1\}$ \\
    $ x_r$ & $=1$ if roof material $r\in\mathcal R$ is selected & $\sum_{r \in \mathcal{R}}  x_r $  $= 1 $, $ x_r $ $\in \{0,1\}$  \\
    $x_w$ & $=1$ if wall material $w\in\mathcal W$ is selected & $\sum_{w \in \mathcal{W}} x_w  = 1, x_w \in \{0,1\}$ \\
    $x^e$ & $=1$ if $e\geq \frac16{{\rm B^\text{fo}}}$ & $x^e\in \{0,1\} $ \\
    $x^{\text{wa}}$ & $=1$ if $l_\text{y}^\text{wa} \geq l_\text{x}^\text{wa} $ & $x^\text{wa}\in \{0,1\} $ \\
      \cmidrule{1-3}
      \multicolumn{3}{l}{\textbf{Changed equations in variable definitions:}}\\
    $F_1$ & seismic load external force & $F_1 = {\rm g}\, {\rm C_\text{d}}v^\text{wa} \sum_{w \in \mathcal{W}} \rho_w x_w$ \\
    $P^\text{fo}_\text{D}$ & self-weight of foundation per unit of lenght & $P^\text{fo}_\text{D} = {\rm g}\, A^\text{fo} \cdot \sum_{f \in \mathcal{F}} x_f\rho_f$ \\
    $P_\text{D}^\text{wa}$ & self-weight of wall per unit length & $P_\text{D}^\text{wa} = {\rm g}\, t^\text{wa} h^\text{wa}\sum_{w \in \mathcal{W}}\rho_w x_w$ \\
    $q^\text{ro}$ & total roof weight & $q^\text{ro} = n^\text{slc}  \cdot v^\text{slc}  \sum_{r \in \mathcal{R}}(x_r\rho_r  + $ \\
    &&\quad ${\rm R^{\text{co}}}\sum_{c \in \mathcal{C}}(x_c \rho_c ))$ \\
      ${\rm \underline{t}^\text{wa}_{w^*}}$ & minimum wall thickness & ${\rm \underline{t}^\text{wa}_{w^*}}= \sum_{w \in \mathcal{W}} {\rm \underline{t}^\text{wa}_{w}} x_w$ \\
         ${\rm \underline{t}^\text{fo}_{f^*}}$ & minimum foundation thickness & ${\rm \underline{t}^\text{fo}_{f^*}} = \sum_{f \in \mathcal{F}} {\rm \underline{t}^\text{fo}_{f}} x_f$ \\
    \bottomrule
\end{tabular}
\end{center}
\end{table}

\subsection{Objective functions}
\label{S:obj-embodied}

Two optimization criteria are considered in the case study:  (i) material costs and (ii) the embodied energy of the building materials. Embodied energy is the energy required to make 1 $kg$ of the material from its ores or feedstocks. The trade-offs between the two optimization criteria are analyzed in Section~\ref{sec:results}.

To compute the costs incurred from a specific building material $m\in\mathcal{M}$, we simply multiply its cost $C_m$ per m$^3$ with the volume of that material.  For example, for a roof material $r\in\mathcal{R}$, the cost is given by $ v^\text{slc}_\text{tot}\cdot  {\rm C}_r$. 
To derive the formulas in Table~\ref{table:variables} for the total volumes we took into account that some components, such as interior walls in the $x$-direction, are shared between multiple rooms.
The material cost for the entire building is the sum of the costs for all materials used:
\begin{equation}\label{eq:cost}
\begin{split}
\text{Cost}\,=\, & v^\text{slc}_\text{tot} \sum_{r \in \mathcal{R}}  {\rm C}_r x_r  + v^\text{co}_\text{tot}\sum_{c \in \mathcal{C}} {\rm C}_c x_c 
+ v^\text{wa}_\text{tot} \sum_{w \in \mathcal{W}} {\rm C}_w x_w \\ &
+ v^\text{fo}_\text{tot} \sum_{f \in \mathcal{F}} {\rm C}_f x_f  +l^\text{re}_\text{tot} {\rm C^\text{re}} .
\end{split}
\end{equation}

Similarly, for a specific material $m\in\mathcal{M}$, the embodied energy per m$^3$ is given by $\rm E_{m}\rho_{m}$, and its overall embodied energy is obtained by multiplying this with the corresponding volume:
\begin{equation}\label{eq:ee}
\begin{split}
\text{Embodied Energy}\, =\,& v^\text{slc}_\text{tot}  \sum_{r \in \mathcal{R}}  {\rm E}_r\rho_rx_r  + v^\text{co}_\text{tot}  \sum_{c \in \mathcal{C}}  {\rm E}_c \rho_c x_c
+ v^\text{wa}_\text{tot}  \sum_{w \in \mathcal{W}} {\rm E}_w \rho_w  x_w\\
& + v^\text{fo}_\text{tot} \sum_{f \in \mathcal{F}} {\rm E}_f \rho_fx_f + l^\text{re}_\text{tot}  {\rm E^\text{re}} \rho^\text{re}.
\end{split}
\end{equation}

Data for unit embodied energy and other sustainability metrics were obtained from ANSYS Granta EduPack software, ANSYS, Inc., 2021 \cite{ansys}. Data on a large number of practical materials were published by Geoff Hammond \& Craig Jones at the University of Bath \citep{hammond2008embodied}. These authors also pair embodied energy values, as $MJ/kg$, with $CO_2$ released in producing the material, as $CO_2/kg$.

\subsection{Formulation of logical conditions}
\label{S:logical}

Some of the parameters and variables in Tables~\ref{table:parameter}--\ref{table:variables} are defined in terms of the building material chosen for the construction.  For example, in the range restriction of the wall thickness $t^\text{wa}$, that is 
\[
{\rm \underline{t}^\text{wa}_{w^*}} \leq t^{\text{wa}} \leq {\rm \bar{t}^\text{wa}},
\]
the minimum ${\rm \underline{t}^\text{wa}_{w^*}}$ depends on the material $w^*\in\mathcal W$ that is chosen for the construction.  To include this condition in the optimization problem, we substitute 
\[
{\rm \underline{t}^\text{wa}_{w^*}} = \sum_{w \in \mathcal{W}} {\rm \underline{t}^\text{wa}_{w}} x_w.
\]
Table~\ref{table:linear} lists all such replacements.

The maximum operator in \eqref{eq:openings1:raw} and \eqref{eq:openings2:raw} cannot be handled by the optimization solver directly and needs to be split into two separate set of inequalities, depending on whether $l_\text{x}^\text{wa}$ is larger than $l_\text{y}^\text{wa}$.
To do this, a binary variable ($x^{\text{wa}}$) is introduced that is 1 when $l_\text{y}^\text{wa}\geq l_\text{x}^\text{wa}$ and 0 otherwise. 
It can easily be verified that the following constraint forces $x^{\text{wa}}$ to take on the correct value:
\begin{equation*}
	0 \leq l_\text{x}^\text{wa} - l_\text{y}^\text{wa} + {\rm M^\text{wa}} \cdot x^{\text{wa}} \leq {\rm M^\text{wa}},
\end{equation*}
where ${\rm M^\text{wa}}={\rm \bar{l}^\text{fl}}-{\rm \underline{l}^\text{fl}}$.
Now the constraints  \eqref{eq:openings1:raw} and \eqref{eq:openings2:raw} can be replaced by the linear constraints
\begin{align*}
w^\text{do} &\leq \tfrac{1}{2} \cdot \big(l_\text{y}^\text{wa} \cdot x^{\text{wa}} + l_\text{x}^\text{wa} \cdot (1-x^{\text{wa}}) \big) \\
l^\text{wi} &\leq \tfrac{1}{2} \cdot \big(l_\text{y}^\text{wa} \cdot x^{\text{wa}} + l_\text{x}^\text{wa} \cdot (1-x^{\text{wa}}) \big).
\end{align*}

In a similar manner, we can express the conditional constraints in Section~\ref{S:foundation_constraints}.  A binary variable $x^e$ is introduced that takes on the value 1 if $e\geq \frac16{{\rm B^\text{fo}}}$ and 0 otherwise.   The following constraint makes sure that $x^e$ indeed satisfies this definition:
\begin{equation*}
0 \leq \frac16{{\rm B^\text{fo}}} - e + \frac16{{\rm B^\text{fo}}}\cdot x^e  \leq \frac16{{\rm B^\text{fo}}}.
\end{equation*}
Now the conditional constraints in Section \ref{S:foundation_constraints} can be expressed as

\begin{equation*}
\begin{split}
(1-x^e) &\cdot \frac{3}{2\, l_\text{x}^\text{wa}}\cdot(F^\text{fo}_\text{D,x}+F_\text{L,x}^\text{wa}) \cdot \left(\frac{1}{{\rm B^\text{fo}}} + \frac{6 \, t^\text{fo} \, e}{({\rm B^\text{fo}})^3} \right) \quad \\
+ \quad x^e &\cdot \left(1 + \frac{t^\text{fo}}{{\rm B^\text{fo}}} \right)\cdot \frac{F^\text{fo}_\text{D,x}+F_\text{L,x}^\text{wa}}{l_\text{x}^\text{wa} ({\rm B^\text{fo}}-2 e)} \quad \leq {\rm \tau^\text{allw}}
\end{split}
\end{equation*}
and
\begin{equation*}
\begin{split}
(1-x^e) &\cdot \frac{3}{2} \cdot (P^\text{fo}_\text{D} + P_\text{D}^\text{wa})\cdot \left(\frac{1}{{\rm B^\text{fo}}} + \frac{6 \, t^\text{fo} \, e}{({\rm B^\text{fo}})^3} \right) \\
  + \quad x^e&\cdot \left(1 + \frac{t^\text{fo}}{\rm B^\text{fo}}\right) \cdot \frac{P^\text{fo}_\text{D} + P_\text{D}^\text{wa}}{{\rm B^\text{fo}} - 2 e}  \quad \leq {\rm \tau^\text{allw}} .
 \end{split}
\end{equation*}
Note that here the inactive term in the inequality is switched off simply by multiplying it with the binary variable or its negation ($1-x^e$).  This is somewhat in contrast to common practice, where often big-M constraints are used to express such a disjunction. 
In our numerical experiments, however, the formulation above resulted in faster solution times compared to a big-M formulation \cite{nemhauser1988integer}.

\section{Case study}
\label{sec:results}
This section showcases how the optimization framework can be utilized to guide decision makers in a variety of circumstances, using the three-room building described in Section \ref{Sec:BuildingModel} as a case study.

\subsection{Nepal-specific case study considerations} 
To demonstrate the range of questions that the proposed optimization methodology can address, this paper incorporates architectural considerations to obtain an optimal solution for minimizing the cost and environmental impact in a multi-room one-story masonry building in Nepal.

The aim of the optimization is to provide a framework for the engineers to work with construction site managers in Nepal to guide infrastructure planning decisions based on model results. The Nepal Department of Urban Development and Building Construction received the first comprehensive training on the BMEG in November 2015 and has shared it with relevant organizations linked to the National Reconstruction Authority. 
The BMEG \cite{Hettiarachchi2016} has been expanded to include 55 materials and provides information on better environmental practices related to design, planning, storage, use, and disposal for common building materials, as well as key environmental costs and benefits. Quantitative information such as embodied energy, CO$_2$ footprint, water usage and several other engineering properties are also provided for each material \cite{BMEG2021}.
The properties of the different building materials which are considered in this paper are given in Table \ref{table:data}. 

The values listed as compressive strength and density are associated with masonry which is composed of two different materials which are joined together layer by layer and side by side. The materials are: (1) the masonry units (the major constituents such as clay bricks, blocks of stone, concrete blocks, pressed earth bricks, etc.); and (2) the mortar paste phase.
For the density, the major constituents (stone, brick, clay brick and concrete/cement block) occupy about 80–95\% \cite{thaickavil2018behaviour} of masonry. Therefore, it was assumed that the density of masonry is 0.875 of unit density of major constituent in masonry plus 0.125 of unit density of mortar paste. For the compressive strength, the relationship $f_p=0.75f_b^{0.75} f_m^{0.31}$ proposed by \cite{lumantarna2014uniaxial} was selected which includes the effects of the properties of the constituents of the masonry and the quality of workmanship. In this equation, $f_p$ is the compressive strength of masonry, $f_b$ is the compressive strength of major constituents in masonry and $f_m$ is the compressive strength of dried mortar paste between masonry units.  

Note that we suppose that different quality grades are available for the wall and foundation materials, where Grade 1 (G1) indicates higher compressive strength and density, and is typically associated with a higher price and greater embodied energy (see Table~\ref{table:data}). The CES EduPack Level 2 Sustainability, Level 3 Sustainability, and Level 3 Eco-Design datasets, technical reports \cite{ashby2018analyzing}, Indian standards and the literature were used as key data sources \cite{ MortarPaste, FiredBrick1, FiredBrick2, CompressedSoil1, CompressedSoil2, UnitofMasonary1, UnitofMasonary2}. Minimum wall and foundation thickness depends on the wall and foundation material. Stone needs a minimum wall and foundation thickness of 0.35 m, while that of brick is	0.23 m. The minimum thickness for concrete and soil blocks is 0.3 m.

Global average material prices were obtained from the CES EduPack \cite{ashby2018analyzing}. Since material prices, as well as properties and environmental impact, may vary significantly with location, it is essential to account for them while making reconstruction decisions. Nepal-specific prices and environmental concerns were obtained from local government officials and industry contacts. In the future, prices in other countries can be obtained by contacting local agencies and government officials involved in reconstruction with help from regional WWF offices. Also, it is notable that during an acute crisis, prices of materials may soar uncontrollably.

We stress that the purpose of this case study is not to give specific recommendation for material selection for this building type. Instead, the following discussion highlights the kind of information that can be obtained with the optimization framework.

\begin{table}[ht]
 \begin{center}
 \begin{minipage}{250pt}
  \caption{Properties of the building materials considered in the case study.}\label{table:data}
 \begin{tabular}{@{}lcrrrr@{}}
  \toprule
Material&Symbol& Density & Cost & EE$^*$ & ${\rm \sigma}^\text{allw,c}_{w}$ \\
&& (kg/m$^3$)& (\$/m$^3$) & (MJ/kg) & (MPa)\\
 \midrule
 \multicolumn{5}{l}{\textbf{Wall and Foundation:}}\\
Stone G1 & (St1) & 2,660 &	35 &	8.71 & 16.07\\
Stone G2 & (St2) & 2,620	& 25 &	7.74 & 10.60\\
Brick G1 & (Br1) & 1,710 & 176 & 12.41& 5.97\\
Brick G2 & (Br2) & 1,560 & 79 & 11.17& 2.33\\
Concrete G1 & (Co1) & 1,730	& 824 & 3.21& 3.26\\
Concrete G2 & (Co2) & 1,210	& 702 & 2.22& 1.88\\
Soil G1	& (So1) & 1,440 &	155 & 0.70& 2.95\\
Soil G2	& (So2) & 1,330 &	145 & 0.67& 2.39\\
 \cmidrule{1-6}
\multicolumn{5}{l}{\textbf{Roof:}}\\
Wood & (Wo) & 850 & 276 & 10.99 &-\\
Bamboo & (Ba) & 1,160 & 803	& 0.15&-\\
 \cmidrule{1-6}
 \multicolumn{5}{l}{\textbf{Roof Cover:}}\\
Plywood & (Pl) & 700 & 385 & 26.83&-\\
Bamboo & (Ba) & 600 & 803 & 0.15 &-\\	
    \bottomrule
  \text{\tiny * Embodied Energy}&&&&&
\end{tabular}
\end{minipage}
\end{center}
\end{table}

\subsection{Multi-objective optimization}
\label{sec:multi-obj}

In Section~\ref{S:obj-embodied}, we defined two optimization criteria: material costs and embodied energy.  Ideally, one would like to minimize both at the same time, but since materials with less embodied energy are often more expensive or have inferior material properties, these objectives are often competing with each other and a decision maker has to strike a balance between the two.

To provide guidance, we are interested in feasible solutions (values of the optimization variables that satisfy all constraints) that are ``Pareto-optimal'', which means that there are no other feasible solutions that provide improvement in both measures at the same time. 
In other words, these solutions represent the best trade-off options since any improvement in one measure would require a compromise in the other.
Our goal is to compute all Pareto solutions and present them as plots, typically referred to as Pareto fronts \cite{berube2009exact}, to the decision maker who can then choose a suitable compromise based on some additional considerations.

We obtained the results in this section using the AMPL optimization modeling language \cite{fourer1990modeling} and the mixed-integer nonlinear programming solver BARON version 21.1.13 \cite{sahinidis:baron:21.1.13} on a M1 Pro with 10-core CPU machine. All the instances are solved to optimality on average in 1.7 seconds. Table~\ref{tab:soltime} in Appendix \ref{app_b} lists the total number of instances solved to generate the different Pareto fronts, together with computation times. The average time per instance is less than 2 seconds in most cases and never exceeds 2.5 seconds. 

\begin{figure}[ht]
\begin{center}
\includegraphics[width=.9\textwidth]{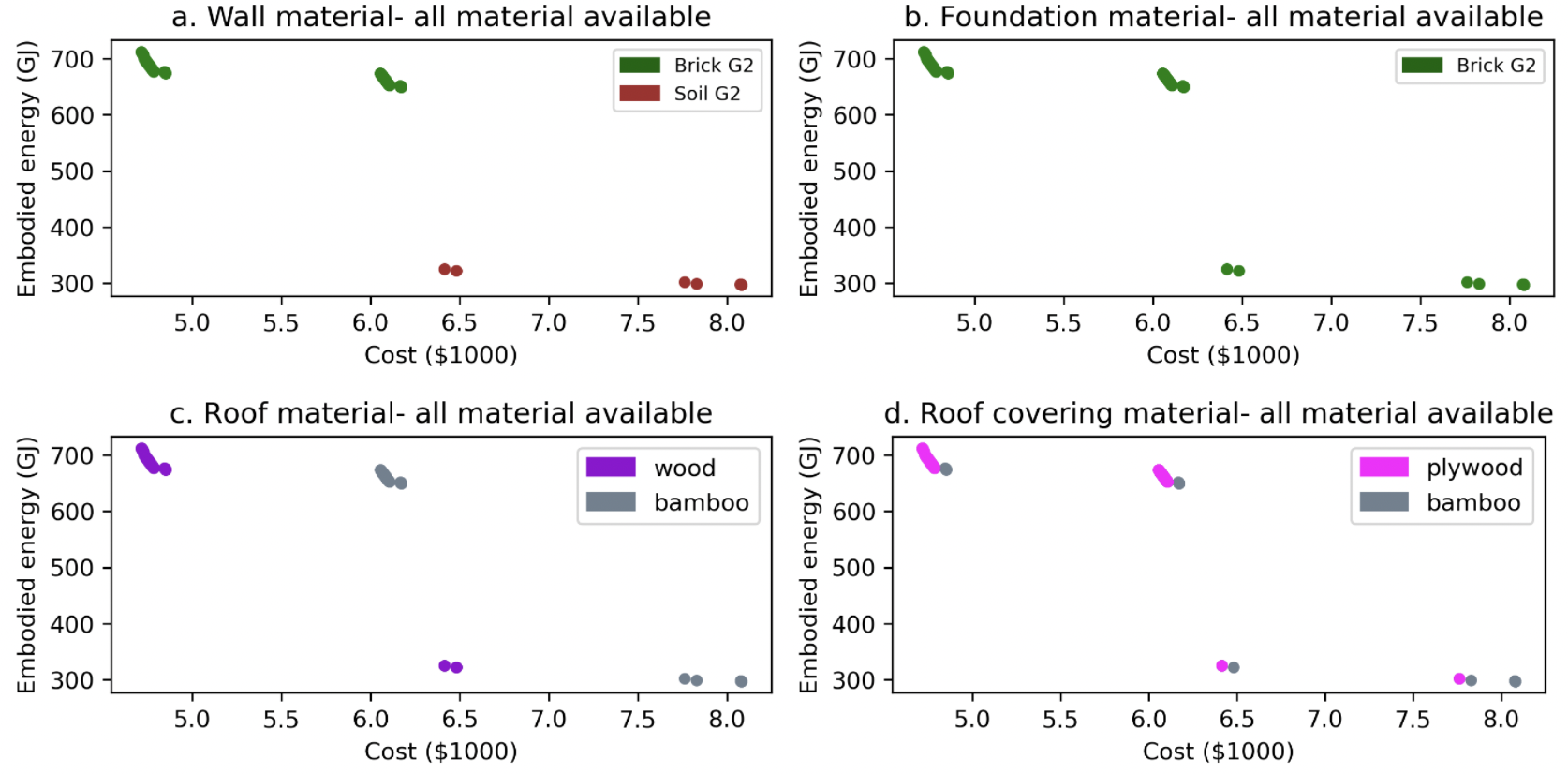}
\end{center}
\caption{
The color-coded Pareto front showing the best budget and embodied energy trade-off. The individual plots depict the selected materials for wall, foundation, roof, and roof cover.}
\label{fig:PF-3room-allMat}
\end{figure}

One way to look at this is to ask the question how the building should be constructed if we want to minimize the embodied energy \eqref{eq:ee} with the restriction that the costs \eqref{eq:cost} have to stay within a given budget limit.
Figure \ref{fig:PF-3room-allMat} has been generated in that fashion, using the $\epsilon$-constraint method \cite{berube2009exact}.
For a wide range of values for the available budget, the optimal solutions were computed and points corresponding to the actual costs (which can be less than the available budget) and embodied energy in the building materials were added to the plots in Figure \ref{fig:PF-3room-allMat}.  The position of the points is identical in the four plots, but their color indicates which building materials were chosen for the different building components. Details of the Pareto optimal solutions including cost, embodied energy, materials, and dimensions are shown in Table \ref{table:summarymap}. 
\begin{table}
\begin{tiny}
 \begin{center}
 \begin{minipage}{350pt}
  \caption{Details of the Pareto-optimal solutions. Each row corresponds to a point or a cluster of points in the Pareto fronts in Figure~\ref{fig:PF-3room-allMat} and gives rise to different design configurations. The precise values for costs and embodied energy are given, as well as the material choices for the different building components.  The last four columns list structural dimensions that distinguish Pareto-optimal solutions within a design. }
  \label{table:summarymap}
 \begin{tabular}{@{}p{0.03\textwidth} p{0.1\textwidth} p{0.1\textwidth}p{0.02\textwidth}p{0.03\textwidth} p{0.02\textwidth} p{0.05\textwidth}p{0.02\textwidth}p{0.08\textwidth} p{0.08\textwidth} p{0.1\textwidth} @{}}
  \toprule
  Design & 
Cost (\$)&
EE$^*$ (GJ) &
Wall &
Found.&
Roof&
Covering&
$n^\text{slc}$ &
$w^\text{do}$ (m)&
$l^\text{wi}$ (m)&
$v^\text{wa}_\text{tot}$ (m$^3$)\\ 
     \midrule
A & 4,715-4,785& 712-677 & Br2 & Br2 & Wo & Pl & 7& 1.10-1.88 & 1.01-1.35 & 23.83-21.77\\ 

B & 4,846-4,852& 677-674 & Br2 & Br2 & Wo & Ba & 7& 1.79-1.88 & 1.35 & 21.92-21.77\\

C & 6,056-6,107& 674-652 & Br2 & Br2 & Ba & Pl & 7& 1.16-1.88 & 1.35 & 22.96-21.70\\ 

D & 6,167-6,173& 652-649 & Br2 & Br2 & Ba & Ba & 7& 1.80-1.88 & 1.35 & 21.85-21.70\\

E & 6,414 &	326 & So2 & Br2 & Wo & Pl	& 7	& 1.89	& 1.35	& 22.80\\
F & 6,481 &	323 & So2 & Br2 & Wo & Ba	& 7	& 1.89	& 1.35	& 22.80\\
G & 7,761 &	302 & So2 & Br2 & Ba & Pl	& 7	& 1.89	& 1.35	& 22.80\\
H & 7,828 &	299 & So2 &	Br2 & Ba & Ba	& 7	& 1.89	& 1.35	& 22.80\\

I & 8,076-8,081 & 298-297 & So2 & Br2 &	Ba & Ba	 & 8 & 1.90-1.89 & 1.35 &	22.56-22.54\\
\bottomrule
 \text{\tiny * Embodied Energy}&&&&&&&&&&
\end{tabular}
\end{minipage}
\end{center}
\end{tiny}
\end{table}

From Figure (4) we see that the minimum budget required to construct the building with any type of material is \$4,715 with Design A, and in that case, the embodied energy cannot be less than 712 GJ.
On the other hand, any construction will need at least 297 GJ of embodied energy, using Design I, and to achieve this, \$8,081 are required.  Having a larger budget available will not make it possible to reduce the embodied energy any further using the materials shown in Table \ref{table:data}.

Second, the best trade-off between the two criteria appears to be achieved by Design E, which corresponds to the left-most point in the cluster around \$6,500 in the Pareto front.  Spending less than \$6,414 will increase the amount of embodied energy by a factor of 2.  The decrease in embodied energy achieved by spending more than \$6,414 is marginal and likely does not justify higher expenses.  The embodied energy of this design is already within 10\% of the smallest possible value.

On the other hand, if \$6,414 are not available, Design D with a cost of \$6,173 appears to be the best choice.  Spending any more money while it is less than \$6,414 would reduce the amount of embodied energy only marginally by less than 10\%.
By reducing the size of the door and the window to their lower limits in Design A, it is possible to save another \$70 since this decreases the amount of expensive rebar material.  This moderate 1.5\% reduction in costs, however, might not justify the resulting 5.2\% increase in embodied energy, and building inhabitants would likely prefer larger wall openings.
In the end, however, it is the judgment of the decision maker to weigh the different trade-offs against each other and choose the final design.

Third, some interesting patterns in the best material selections can be observed:  (i) whenever bricks are used, the cheapest grade is sufficient to meet the structural requirements; (ii) the key factor in reducing the embodied energy is the switch from bricks to soil blocks for the wall material when it can be afforded; (iii) wood is the preferred roof material since bamboo is significantly more expensive but reduces the embodied energy only marginally, and (iv) for the same reason, plywood is the preferred roof covering material.

\subsection{Unavailability of some materials}
\label{Unavailabilitymat}
\begin{figure*}[ht]
\begin{center}
\includegraphics[width=1\textwidth]{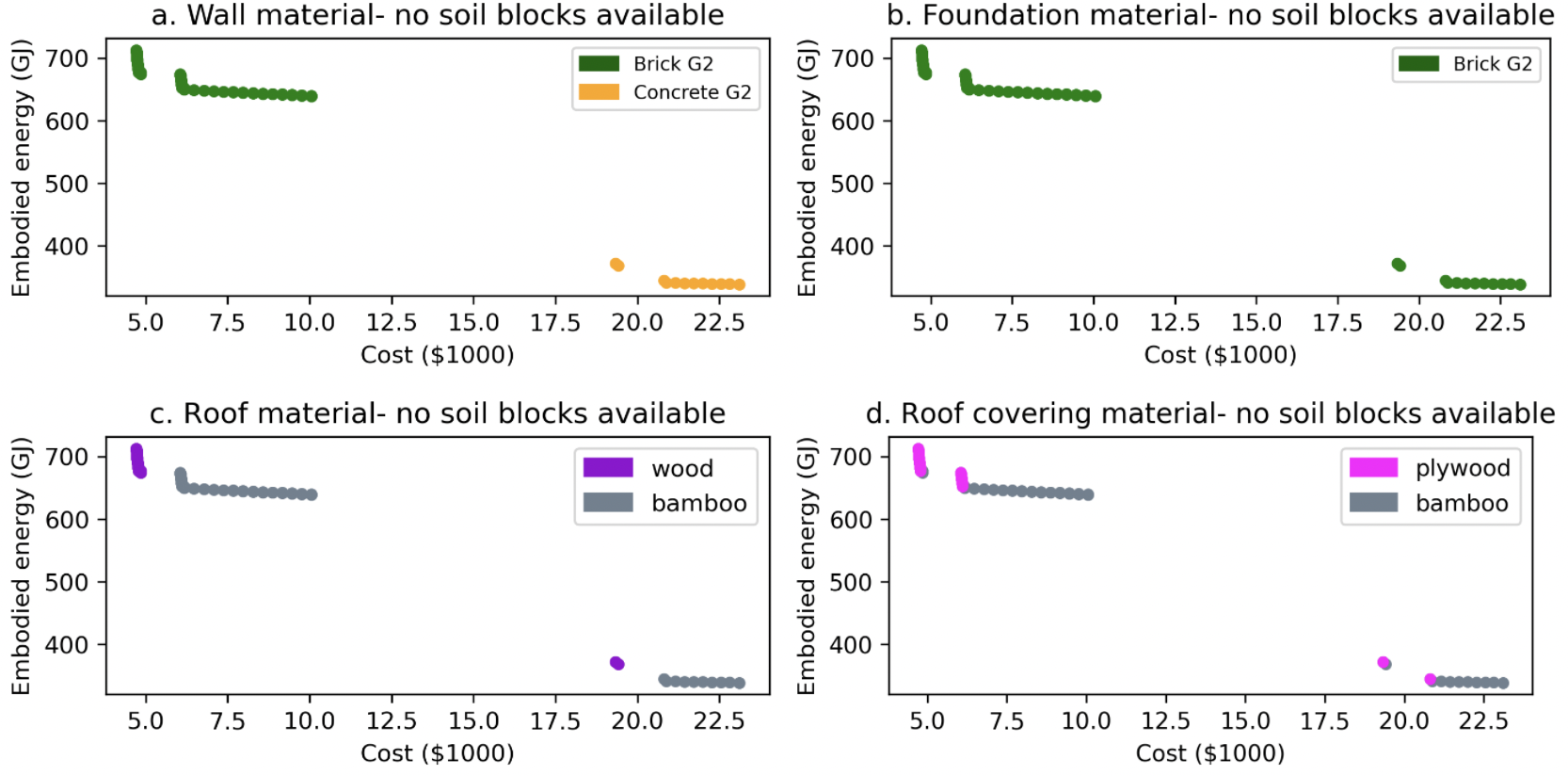}
\end{center}
\caption{ The color-coded Pareto front showing the best budget and embodied energy trade-off when soil block is not available as building material. The individual plots depict the selected materials for wall, foundation, roof, and roof cover.}
\label{fig:MaterialSelection-3-room-exceptsoil}
\end{figure*}

Factors such as climate change and transportation network disruption due to disasters are likely to affect material availability and price in the short and long term. First, regional changes in precipitation and drought patterns \cite{rodell2018emerging,diffenbaugh2017quantifying,pachauri2014climate} will likely impact the availability of water-dependent resources including timber, bamboo, cement blocks and brick. Second, the growing frequency and intensity of storms will likely inflict more damage on regional forestry resources. Third, a changing climate is likely to impact forest health by expanding the likelihood of damaging pest and disease outbreaks \cite{kirilenko2007climate}, which could greatly alter the types of timber available for construction.
The discussion in the previous section revealed that the possibility of using soil block as the wall material is important to reduce the amount of embodied energy.
A natural question to ask is what should be done when soil block is not available.
Figure \ref{fig:MaterialSelection-3-room-exceptsoil} shows the Pareto front for this setting.  A significant decrease in the amount of embodied energy can still be achieved by using concrete block instead of bricks for the wall material, but at a very high cost.  It requires \$19,329 to reduce the amount of embodied energy to less than 400 GJ, and even with an unlimited budget, the amount of embodied energy cannot be decreased below 339 GJ, which is about 14\% worse than what could be obtained when soil block is used.
The other conclusions from the previous section are still valid:  The cheapest brick grade should always be chosen as foundation material, and wood and plywood are preferable choices for roof and roof covering, respectively.

\begin{figure*}[ht]
\begin{center}
\includegraphics[width=1\textwidth]{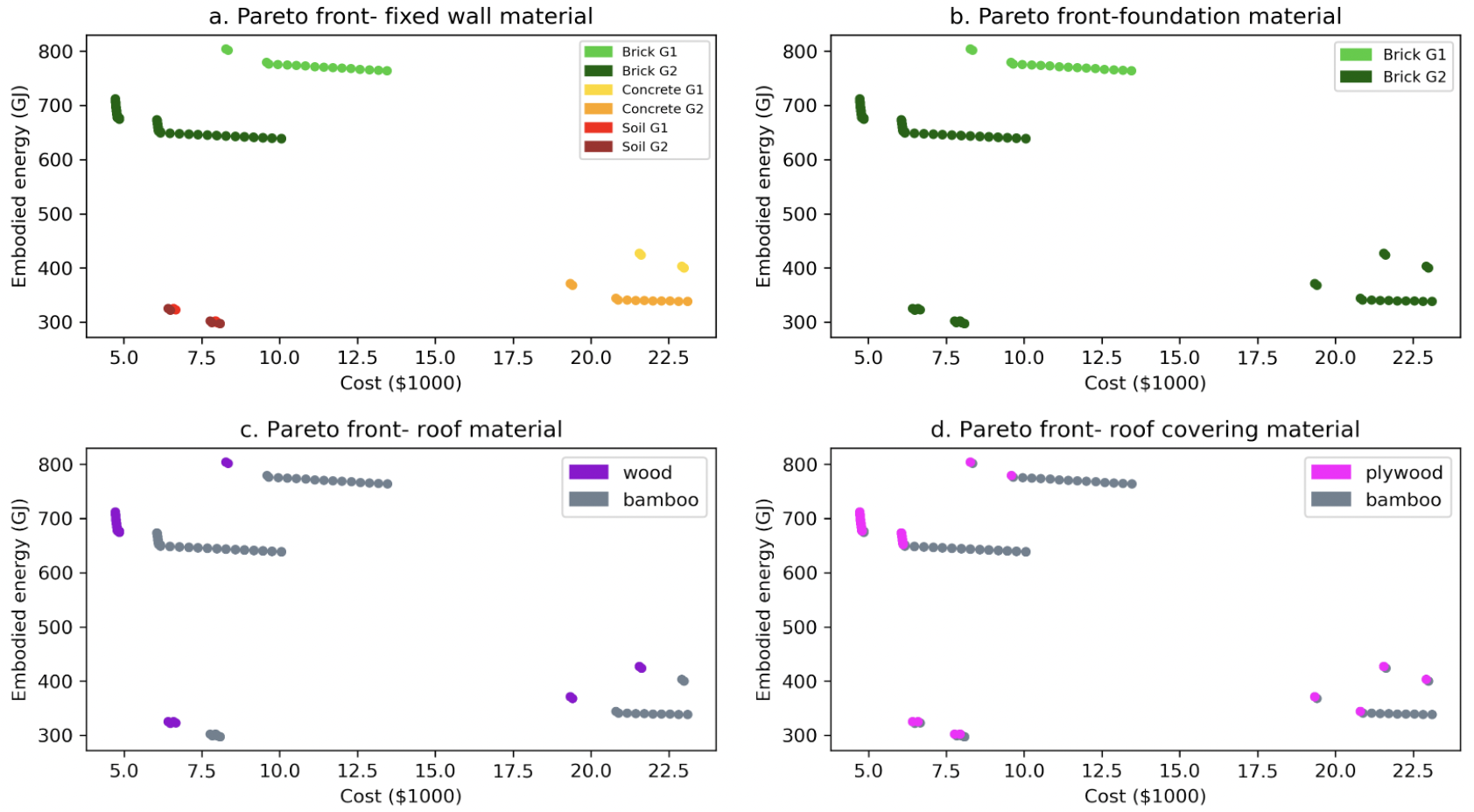}
\end{center}
\caption{ Combined Pareto fronts with fixed wall material.}
\label{fig:pareto-wallmat-allinall}
\end{figure*}

To better understand the effect of the different options for the wall material in general, we combined several Pareto fronts in Figure \ref{fig:pareto-wallmat-allinall}.  The different colors in Figure \ref{fig:pareto-wallmat-allinall}.a represent the Pareto optimal solutions when the wall material is fixed to a particular choice.  For instance, the dark green points labeled Brick G2 correspond to the configurations that provide the lowest embodied energy for a given budget that can be achieved when Brick G2 is used for the wall.  As before, the colors in the remaining three plots identify the best material selection for the other building components.
For consistency, we only permit bricks of the same grade for the foundation and the walls.

Figure \ref{fig:pareto-wallmat-allinall} shows that the worst material in terms of embodied energy is brick G1. 
If soil block is not available, concrete block can be used to reduce the amount of embodied energy, but with a cost that is up to 2.5 times that of the other materials.
It turns out that stone cannot be used for the wall material, although it is the cheapest option. The allowable compressive strength and density properties of the available stone grades cannot meet the building structure requirements. 
We explore this observation in more detail in Section~\ref{sec:thresh}.

\subsection{Sensitivity to prices}
Prices are the most uncertain parameters in these calculations, particularly during an acute crisis, and it is important to understand how price changes impact the optimal design decision.
Economies can be devastated by disasters, causing dramatic distortions in local prices and material availability \cite{rose2017economic}. Prices and material properties will vary in different regions and countries.
As an example, we look at the prices of wall materials which make up the largest fraction of the overall costs.

Consider a price change for soil block G2.  The change only affects the overall construction costs \eqref{eq:cost} of the different designs that involve soil block G2.  The designs themselves do not change and will still satisfy all engineering and safety restrictions, and the amount of embodied energy remains the same.  If ${\rm C}^\text{orig}_{w^*}$ 
is the original and ${\rm C}^\text{new}_{w^*}$ is the new price for soil block G2, the cost of a design involving soil block G2, such as Design E, can be recalculated from \eqref{eq:cost} as
\[
\text{Cost}_F^\text{new} = \text{Cost}_F^\text{orig} + v^\text{wa}_\text{tot} ({\rm C }^\text{new}_{w^*} - {\rm C}^\text{orig}_{w^*}).
\]
Graphically, in Figures~\ref{fig:PF-3room-allMat} and \ref{fig:pareto-wallmat-allinall}, an increase in price shifts the red points for soil blocks G2 horizontally to the right.  
This insight permits us to update the Pareto front in Figure~\ref{fig:PF-3room-allMat} using the combined Pareto fronts in Figure~\ref{fig:pareto-wallmat-allinall}:  After shifting the red points horizontally according to the price change, we choose all non-dominated points in Figure~\ref{fig:pareto-wallmat-allinall}; i.e., the points for which there is no other point that is better with regards to both criteria.

For instance, suppose that a maximum budget of \$8,000 is available.  Then Design E, with original cost \$6,414 and total wall volume of 22.80 $m^3$ is still a good choice as long as the price is not larger than
\begin{align*}
{\rm C}^\text{new}_{w^*} = \;& {\rm C}^\text{orig}_{w^*} + \frac{\text{Cost}_F^\text{new}-\text{Cost}_F^\text{orig}} {v^\text{wa}_\text{tot}} = \\ & 145 +
\frac{8,000-  6,414} {22.80}
= \$214.56.
\end{align*}
Once the price exceeds this threshold, only designs using brick can be paid for.

\subsection{Varying the floor size}

The analysis so far has considered buildings in which each room has a floor area of $A^{\text{fl}}_{\text{min}}= \text{10 m}^2$.  A decision maker might also want to take into account flexibility with respect to the room size as a third optimization criteria: besides minimizing costs and embodied energy, we also want to maximize the size of the room.  Now there are three objectives that are competing with each other.

\begin{figure}[ht]
\begin{center}
\includegraphics[width=.55\textwidth]{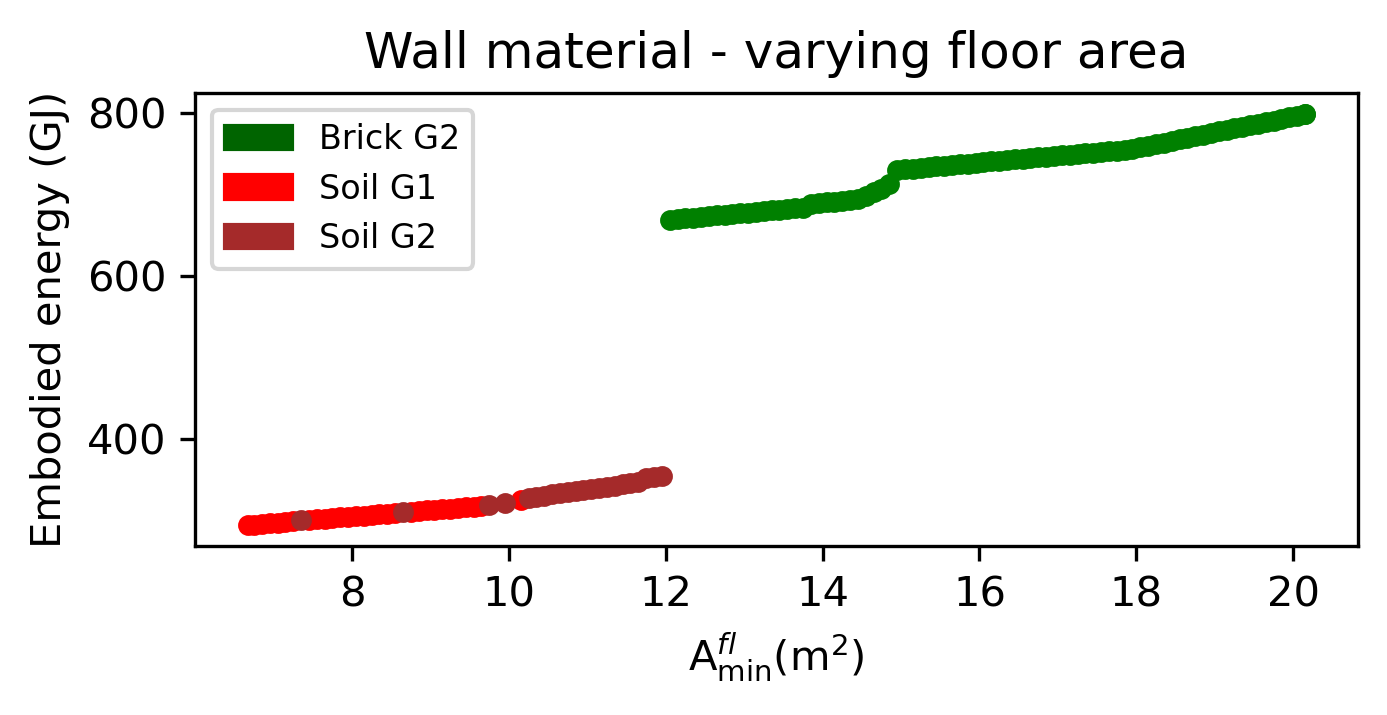}
\end{center}
\caption{ The color-coded Pareto front showing the best floor area and embodied energy trade-off for a budget fixed at \$7,000. It depicts the selected materials for wall.}
\label{fig:PF-area}
\end{figure}

Figure~\ref{fig:PF-area} shows the Pareto curve for embodied energy vs.\ floor area (instead of cost) when the budget is fixed to \$7,000. 
As before, the Pareto front is obtained with the $\epsilon$-constraint method \cite{berube2009exact} by minimizing the embodied energy while varying the floor area.

We see that soil blocks remains the preferred wall material up to a floor area of 11.95 m$^2$. For larger buildings, however, the use of soil blocks would exceed the available budget and Brick G2 becomes the best option, at the expense of a considerable increase of 316 GJ in the embodied energy.

For the floor area less than 10.2 m$^2$, soil block G1 is preferred except in some cases where using soil block G2 is a better choice. This choice of soil block G2 results in higher volume of superstructure and foundation due to the increase in the thickness of wall in the corresponding designs. However, it requires slightly less total budget due to the lower price of soil block G2. 

It is important to note that a stakeholder can recalculate the Pareto front for higher or lower budgets as well.  Furthermore, one could study the trade-off between floor area and budget with a given fixed limit on the embodied energy.

In general, the multi-objective optimization approach makes it possible to examine the trade-off between a variety of criteria at the same time.

\subsection{Understanding the impact of safety thresholds}
\label{sec:thresh}

As discussed in Section~\ref{Unavailabilitymat}, no design that satisfies the safety constraints in Section~\ref{S:foundation_constraints} can include stone as the wall material.
The threshold parameters in these constraints are based on worst-case values that cover a wide geographical area.
Here we explore whether stone, which is much cheaper than all other materials, becomes a viable option if one or more of these thresholds is relaxed.

One crucial parameter in the safety constraints is the width of the foundation (${\rm B^\text{fo}}$), for which the conservative value 0.8$m$ had been chosen, based on the worst possible soil characteristics in the region.
The proposed optimization framework can be used to determine a threshold for ${\rm B^\text{fo}}$ above which stone is a viable wall material.  
This is done by minimizing ${\rm B^\text{fo}}$ with the choice of wall material fixed to either grade of stone, subject to all engineering and safety constraints, but without restrictions on costs and embodied energy.
These calculations find that both stone G1 and G2 are feasible materials in geographical regions where ${\rm B^\text{fo}}$ is at least 0.81$m$.

\begin{figure*}[ht]
\begin{center}
\includegraphics[width=1\textwidth]{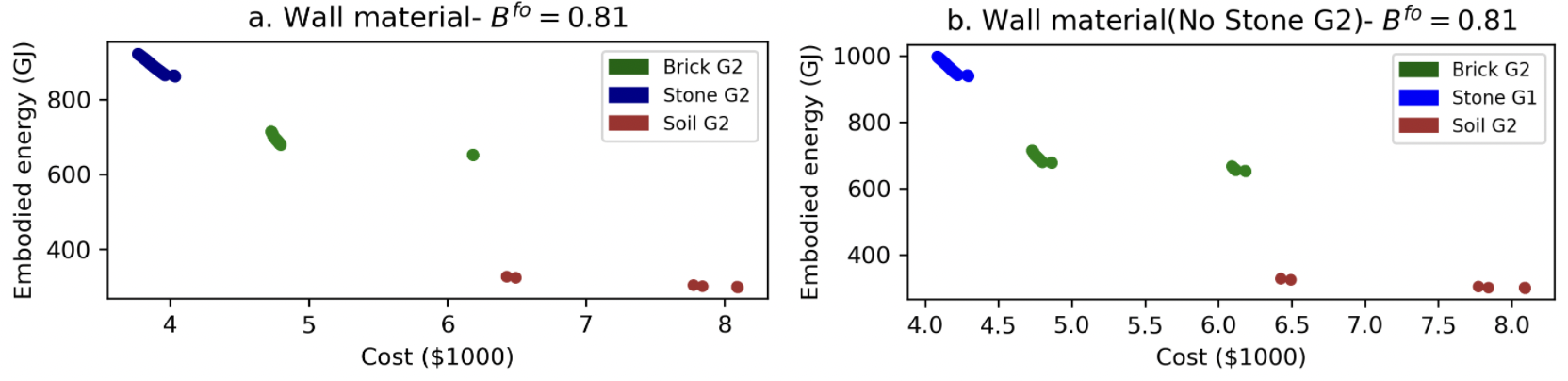}
\end{center}
\caption{ The color-coded Pareto front showing the best budget and embodied energy trade-off when foundation width is 0.81 rather than 0.8. It depict the selected materials for wall. Plot a: All materials are available. Plot b: Stone G2 is not available.}
\label{fig:PF-0.81}
\end{figure*}

Figure \ref{fig:PF-0.81} shows the Pareto fronts, updated for ${\rm B^\text{fo}}=$ 0.81$m$. In Figure \ref{fig:PF-0.81}a, all materials are available. 
In contrast to the original Pareto front in Figure~\ref{fig:PF-3room-allMat}, we now see designs that use Stone G2, with costs ranging from \$3,771 to \$4,035 and embodied energy between 922 and 862 GJ.

Figure \ref{fig:PF-0.81}b depicts the result for the case when stone G2 is not available. It shows an increase in total cost and embodied energy when stone G1 is used.
The embodied energy varies between 997 and 939 GJ, and the cost ranges from \$4,084 to \$4,293.

This study reveals that a small change in the foundation width safety parameter opens the door to use building materials that were not permitted when strict worst-case thresholds are used.  This can be a very important insight, since it might lead to different recommendations that advocate the potential use of stone in most of the geographical regions, with the exception of some small areas. However, the lower cost would come at the expense of much greater embodied energy. Similar sensitivity studies can be performed for other parameters in the optimization model.

The different scenarios of this case study illustrate the wide range of insights that the proposed framework can provide related to the optimization of building designs under different conditions.

\section{Discussion}
\label{S:7}

This work demonstrates a systematic approach for understanding trade-offs and interpreting the varying solutions that are obtained after changing model parameters and constraints. A summary table such as that shown in Table  \ref{table:summarymap} clearly demonstrates the trade-offs between alternatives. It provides details for the different sections of the Pareto front and the corresponding design configurations. Our approach is distinguished from the existing ones in that it allows rapid screening of alternatives. It facilitates the process for selection of the final building design that we introduced in Figure \ref{fig:flowchart}
and represents a crucial improvement. That is because, as Lesken et al. \cite{leskens2014decisions} mentioned, “in the disaster recovery phase, problems do not emerge with well-defined properties, and after a disaster, effective alternatives are vague.” 

The inverse relationship generally found between the cost of a construction project and the project’s embodied energy of materials suggests that site managers must work with different stakeholders to optimize the use of available resources and meet demands, while still taking sustainability into consideration. 
However, the majority of previous works on post-disaster management focused either on minimizing the cost or the required service and reconstruction time. Zamanifar and Hartmann \cite{zamanifar2020optimization} performed an extensive review of optimization-based decision-making models for disaster recovery planning until 2020. According to their analysis, traffic flow, traffic time and economic factors were the main focus of these models. Out of 241 papers that they fully investigated, there was no work on reconstruction at the building level. As an example, Majumder et al.\ \cite{majumder2020review} developed an optimization problem to minimize debris collection cost without considering the debris to be used for reconstruction. 
One paper by Onan et al.\ \cite{onan2015evolutionary} that considered both environmental impacts and costs for disaster management optimized transportation cost and hazardous waste exposure risk. However, their approach does not include rapid screening of alternatives to account for the situations where some scenarios are not available. 

Another important contribution of this paper is that it provides optimal and sustainable detailed building designs. A systematic review paper on temporary housing management modeling and decision-making methods by Perrucci and Baroud  \cite{perrucci2020review} showed that expenditure, displacement duration, and health and social well-being were the three main categories evaluated by researchers. Also, according to Perrucci’s and Baroud’s analysis, most of the optimization papers have focused on supply chain and logistics management, rather than the detailed building structure and material selection. El-Anwar et al.\ \cite{el2016efficient} optimized traffic disruption and reconstruction costs at the same time, but they did not optimize reconstruction costs based on the building design. 
Castro-Lacouture et al.’s \cite{castro2009optimization} model and approach for selecting materials using a LEED-based rating system in Colombia are closely related to our work. However, they relied on LEED credits for material selection and did not determine the building design parameters.

Another contribution of our paper in the post-disaster management literature is proposing a quantitative optimization model that considers several variables and makes several optimal decisions at the same time. The scheme is versatile, in that changes to any, or even all, of the variables can be taken into account, which makes it different from existing works. For example, in Majumder et al.\ \cite{majumder2020review} the decision variable is the quantity of debris that shifts from disaster-affected regions to debris management sites. Also, in Onan et al.’s paper \cite{onan2015evolutionary} the decision variables are temporary storage location and transportation of the waste in a zone to a temporary storage site. Their mathematical optimization model is linear integer and can be solved in reasonable time.
Other papers that consider a greater number of decision variables use non-exact methods to solve their problem. As an example, Ghannad et al.\ \cite{ghannad2020multiobjective} integrated the analytical hierarchy process (AHP) and the multi-objective genetic algorithm approach for a comprehensive consideration of qualitative and quantitative factors. However, their model lacks sustainability considerations and they found approximate, rather than exact, solutions by using a heuristic method. 

It is also notable that, in contrast with some game theoretic modelling approaches where different stakeholders are modelled without their involvement and feedback in decision-making (e.g., \cite{yue2014game}), direct stakeholder engagement is considered in the process of decision-making with our approach (Figure \ref{fig:flowchart}). This paper’s framework is designed to be used by people in the field and real stakeholders. A major frontier in sustainable mathematical optimization R\&D is bringing communities into the decision-making process and integrating it as a core concept in the optimization framework, as we have shown in this paper. Roth et al.\ \cite{roth2021importing} bring a similar approach to the process systems engineering community. 

The model presented in this paper can be adapted to different regions in which construction will take place. Different regions present different availability of materials and, as our results demonstrate, the site soil type can limit structural dimensions. While many parameters used in this model are based on worst-case scenarios across Nepal, we observed that a small increase in soil stability allowed for the use of stone for the wall material that greatly reduced the cost. Thus, site managers might use this knowledge to guide stakeholders for construction to perform an in-depth site soil analysis to make construction less expensive and/or more environmentally friendly.

Future advances to improve the building model include the integration of ArcGIS \cite{law2015getting} to obtain soil data from a multitude of regions and adapting the model to find all possible feasible solutions for construction in these regions and the trade-offs that come with them. Partnering with site construction managers for a case study trial of the model presented here would be crucial for the model to pass a proof of concept and real-world usability. Not only would site construction access allow for more complete inclusion of site data, but it would allow for an understanding of the timeline and decision-making process for construction from which we can begin to understand at what point(s) model results would be most useful for planning. 

Future analyses can consider metrics beyond embodied energy including life-cycle greenhouse gas emissions and water consumption to provide additional insight. Also, modeling qualitative data and demands surrounding the stakeholders of a given construction project as part of a multi-objective optimization would also serve as a significant extension of this context. This type of optimization could ideally take into consideration ethical contexts under which labor is obtained or the degree to which one might value sustainability relative to cost.

\backmatter

\bmhead{Acknowledgments}
Vidushi Dwivedi was supported in part by the Resnick Family Social Impact Program of Northwestern's Institute for Sustainability and Energy (ISEN) and the World Wildlife Fund (WWF). Daniel J.\ Garcia was supported by a Northwestern University Presidential Fellowship. Additional support was provided by Leslie and Mac McQuown. 

We thank Mike McMahon for help with coordination and comments on the paper and Jennifer B.\ Dunn for discussions providing a life cycle analysis perspective.

We also thank Ah Kim, Can Divitoglu, and Giovanna Varalta, who created an initial version of the AMPL model as part of their project in the Northwestern course IEMS 394 \emph{IE Client Project Challenge.} 

\begin{appendices}
\label{appendix}

\section{Complete mathematical optimization model}\label{secmodel}

For completeness, the full optimization model is summarized in this section. It uses the symbols and variables given in Tables \ref{table:parameter}, \ref{table:variables} and \ref{table:linear}.

The final objective function that is proposed in Section \ref{S:obj-embodied} is in \eqref{eq:obj:smooth}. Equation \eqref{eq:financ:smooth} is the second objective function that is moved to the constraints by $\epsilon$-constraint method. The full building model constraints that are proposed in Section \ref{Sec:BuildingModel} are demonstrated in \eqref{eq:roof:smooth}-\eqref{eq:opening:smooth}. The equations in \eqref{eq:bound1:smooth}-\eqref{eq:boundn:smooth} are the variable bounds and definitions that are given in Table \ref{table:variables} and \ref{table:linear}.

\setlength{\belowdisplayskip}{1pt} \setlength{\belowdisplayshortskip}{1pt}
\setlength{\abovedisplayskip}{1pt} \setlength{\abovedisplayshortskip}{1pt}
\begin{small}
\begin{flalign}
\label{eq:obj:smooth}
&\min \quad v^\text{slc}_\text{tot} ( \sum_{r \in \mathcal{R}}  x_r {\rm E}_r\rho_r )  + v^\text{co}_\text{tot} ( \sum_{c \in \mathcal{C}} x_c {\rm E}_c \rho_c )\notag\\
&+ v^\text{wa}_\text{tot} ( \sum_{w \in \mathcal{W}} x_w {\rm E}_w \rho_w )+ l^\text{re}_\text{tot}( {\rm E^\text{re}} \rho^\text{re})
 + v^\text{fo}_\text{tot} ( \sum_{f \in \mathcal{F}} x_f {\rm E}_f \rho_f )
\end{flalign}
\end{small}
Subject to,
\\
\textbf{Roof constraints}
\begin{align}\label{eq:roof:smooth}
n^\text{slc} {\rm w^\text{be}} \leq l_\text{y}^\text{wa} \leq  (n^\text{slc}-1) \cdot {\rm \bar{s}^\text{be}} + n^\text{slc} {\rm w^\text{be}}
\end{align}
\textbf{Wall constraints}
\begin{align}
l_\text{x}^\text{fl} \cdot l_\text{y}^\text{fl} \geq {\rm \underline{A}^\text{fl}}
\end{align}
\begin{align}
\frac{F_\text{D,x}^\text{wa}+F_\text{L,x}^\text{wa}}{A^\text{wa}_\text{x}}+\frac{M_{w,\text{x}}}{S_{\text{x}}} \leq  \sum_{w \in \mathcal{W}} {\rm \sigma}^\text{allw,c}_{w} x_w
\end{align}
 \begin{align}
\frac{F_\text{D,x}^\text{wa}+F_\text{L,x}^\text{wa}}{A^\text{wa}_\text{x}}+\frac{M_\text{e,x}}{S_{\text{x}}} \leq \sum_{w \in \mathcal{W}} {\rm \sigma}^\text{allw,c}_{w} x_w
\end{align}
\begin{align}
-\frac{F_\text{D,x}^\text{wa}+F_\text{L,x}^\text{wa}}{A^\text{wa}_\text{x}}+\frac{M_{w,\text{x}}}{S_{\text{x}}} \leq  {\rm \sigma^\text{allw}_{t}}
\end{align}
 \begin{align}
 -\frac{F_\text{D,x}^\text{wa}+F_\text{L,x}^\text{wa}}{A^\text{wa}_\text{x}}+\frac{M_\text{e,x}}{S_{\text{x}}} \leq {\rm \sigma^\text{allw}_{t}} 
\end{align}
\begin{align}
\frac{3 F_{w,\text{x}}}{2 A^\text{wa}_\text{x}}\leq {\rm \tau^\text{allw}}
\end{align}
  \begin{align}
 \frac{3 F_\text{e}}{2 A^\text{wa}_\text{x}}\leq {\rm \tau^\text{allw}}  
\end{align}
\begin{align}
 \frac{P_\text{D}^\text{wa}}{t^\text{wa}}+\frac{M_{w,\text{y}}}{S_{\text{y}}} \leq \sum_{w \in \mathcal{W}} {\rm \sigma}^\text{allw,c}_{w} x_w  
  \end{align}
   \begin{align}
 \frac{P_\text{D}^\text{wa}}{t^\text{wa}}+\frac{M_\text{e,x}}{S_{\text{y}}} \leq \sum_{w \in \mathcal{W}} {\rm \sigma}^\text{allw,c}_{w} x_w 
\end{align}
  \begin{align}
 -\frac{P_\text{D}^\text{wa}}{t^\text{wa}}+\frac{M_{w,\text{y}}}{S_{\text{y}}} \leq {\rm \sigma^\text{allw}_{t}} 
  \end{align}
 \begin{align}
 -\frac{P_\text{D}^\text{wa}}{t^\text{wa}}+\frac{M_\text{e,x}}{S_{\text{y}}} \leq {\rm \sigma^\text{allw}_{t}}  
\end{align}
\begin{align}
 \frac{3 F_{w,\text{y}}}{ 2 A^\text{wa}_\text{y}}\leq {\rm \tau^\text{allw}} 
  \end{align}
 \begin{align}
 \frac{3F_\text{e}}{2 A^\text{wa}_\text{y}}\leq {\rm \tau^\text{allw}}  
\end{align}
\textbf{Foundation constraints}
\begin{align}\label{eq:27}
0 \leq \frac16{{\rm B^\text{fo}}} - e + \frac16{{\rm B^\text{fo}}}\cdot x^e  \leq \frac16{{\rm B^\text{fo}}}
\end{align}
\begin{align}\label{eq:28}
\begin{aligned}
&(1-x^e) \cdot \frac{3}{2 \cdot l_\text{x}^\text{wa}}\cdot(F^\text{fo}_\text{D,x}+F_\text{L,x}^\text{wa}) \cdot \left(\frac{1}{{\rm B^\text{fo}}} + \frac{6 \cdot t^\text{fo} \cdot e}{({\rm B^\text{fo}})^3} \right) \quad \\
&+ \quad  x^e \cdot (1 + \frac{t^\text{fo}}{{\rm B^\text{fo}}} )\cdot \frac{F^\text{fo}_\text{D,x}+F_\text{L,x}^\text{wa}}{l_\text{x}^\text{wa} ({\rm B^\text{fo}}-2 e)} \quad \leq {\rm \tau^\text{allw}} 
\end{aligned}
\end{align}
\begin{align}\label{eq:29}
&(1-x^e) \cdot \frac{3}{2} \cdot (P^\text{fo}_\text{D} + P_\text{D}^\text{wa})\cdot \left(\frac{1}{\rm B^\text{fo}} + \frac{6 \cdot t^\text{fo} \cdot e}{({\rm B^\text{fo}})^3} \right) \notag\\
&+ \quad x^e\cdot (1 + \frac{t^\text{fo}}{{\rm B^\text{fo}}}) \cdot \frac{P^\text{fo}_\text{D} + P_\text{D}^\text{wa}}{{\rm B^\text{fo}} - 2 e}  \quad \leq {\rm \tau^\text{allw}} 
\end{align}
\begin{equation}\label{eq:30}
	0 \leq l_\text{x}^\text{wa} - l_\text{y}^\text{wa} + {\rm M^\text{wa}} \cdot x^{\text{wa}} \leq {\rm M^\text{wa}}
\end{equation}
\textbf{Opening constraints}
\begin{equation}\label{eq:31}
w^\text{do} \leq \frac{1}{2} \cdot (l_\text{y}^\text{wa} \cdot x^{\text{wa}} + l_\text{x}^\text{wa} \cdot (1-x^{\text{wa}}) ) 
\end{equation}
\begin{equation}\label{eq:31_1}
l^\text{wi} \leq \frac{1}{2} \cdot (l_\text{y}^\text{wa} \cdot x^{\text{wa}} + l_\text{x}^\text{wa} \cdot (1-x^{\text{wa}}) ) 
\end{equation}
\begin{equation}
	l^\text{wi} \leq \frac{h^\text{wa}}{2} 
\end{equation}
\begin{equation}\label{eq:opening:smooth}
n^\text{re} {\rm d^\text{re}} + (n^\text{re} - 1) {\rm \underline{s}^\text{re}} \leq t^\text{wa} \\
\end{equation}
\textbf{Financial constraint}
\begin{small}
\begin{equation}\label{eq:financ:smooth}
\begin{split}
v^\text{slc}_\text{tot}  \sum_{r \in \mathcal{R}}  {\rm C}_r x_r   + v^\text{co}_\text{tot} \sum_{c \in \mathcal{C}} {\rm C}_c x_c  
+  
v^\text{wa}_\text{tot}  \sum_{w \in \mathcal{W}} {\rm C}_w x_w & \\+l^\text{re}_\text{tot} {\rm C^\text{re}} 
+ v^\text{fo}_\text{tot}  \sum_{f \in \mathcal{F}} {\rm C_f} x_f & \leq {\rm B^{\text{avail}}}
\end{split}
\end{equation}
\textbf{Variables bounds and equations}
\begin{equation}\label{eq:bound1:smooth}
\{l_\text{x}^\text{wa},l_\text{y}^\text{wa}\} = \{  l_\text{x}^\text{fl}  + 2 t^\text{wa}, l_\text{y}^\text{fl}  + 2 t^\text{wa}\}
\end{equation}
\setlength{\belowdisplayskip}{2pt} \setlength{\belowdisplayshortskip}{2pt}
\setlength{\abovedisplayskip}{2pt} \setlength{\abovedisplayshortskip}{2pt}
\begin{equation}
l^{\text{re}} = n^\text{re} \cdot 2 ( w^\text{do} + {\rm h^\text{do}} ) +4 l^\text{wi}
\end{equation}
\begin{equation}
F_\text{D,x}^\text{wa} = P_\text{D}^\text{wa} l_\text{x}^\text{wa} + (0.5)({\rm g})q^\text{ro}
\end{equation}
\begin{equation}
F_\text{L,x}^\text{wa} = 0.5 {\rm P_\text{L}} l_\text{x}^\text{wa} t^\text{wa}
\end{equation}
\begin{equation}
\{F_{w,\text{x}},F_{w,\text{y}}\} = 0.5 {\rm C_\text{f}} \cdot \{l_\text{y}^\text{wa},l_\text{x}^\text{wa}\} \cdot h^\text{wa} {\rm P^\text{design}}
\end{equation}
\begin{equation}
\{M_{w,\text{x}},M_{w,\text{y}}\} = \{F_{w,\text{x}}\cdot h^\text{wa}, F_{w,\text{y}}\cdot h^\text{wa}\}
\end{equation}
\begin{equation}
v^\text{wa} = t^\text{wa} \cdot (2 \cdot h^\text{wa} \cdot (l_\text{y}^\text{fl}+l_\text{x}^\text{wa}) - (w_d \cdot {\rm h^\text{do}} + (l^{\text{wi}})^2) )
\end{equation}
\begin{equation}
\begin{split}
v^\text{wa}_\text{tot} = {\rm n^{\text{rm}}} v^\text{wa} - ({\rm n^{\text{rm}}}-1) t^\text{wa} h^\text{wa} l_\text{x}^\text{wa}
\end{split}
\end{equation}
\begin{equation}
F_\text{e} = \frac{F_1}{2}
\end{equation}
\begin{equation}
M_\text{e,x} = F_\text{e}\cdot h^\text{wa}
\end{equation}
\begin{equation}
\{A^\text{wa}_\text{x},A^\text{wa}_\text{y}\}=\{t^\text{wa}\cdot l_\text{x}^\text{wa},t^\text{wa}\cdot l_\text{y}^\text{wa}\}
\end{equation}
\begin{equation}
\{S_{\text{x}},S_{\text{y}}\}= \{ \frac16 t^\text{wa} \cdot  (l_\text{x}^\text{wa})^2,\frac16 t^\text{wa} \cdot  (l_\text{y}^\text{wa})^2\}
\end{equation}
\begin{equation}
v^\text{slc} = {\rm A^\text{be}} l_\text{x}^\text{wa} + 2 {\rm A^\text{ra}} \sqrt{ (\frac{{\rm R^{\text{be}}}\cdot l_\text{x}^\text{wa}}{2} )^2 + (\frac{l_\text{x}^\text{wa}}{2})^2}
\end{equation}
\begin{equation}
v^\text{slc}_\text{tot} = {\rm n^{\text{rm}}} \cdot n^\text{slc} \cdot  v^\text{slc}
\end{equation}
\begin{equation}
v^\text{co}_\text{tot}= {\rm R^{\text{co}}}\cdot v^\text{slc}_\text{tot}
\end{equation}
\begin{equation}
A^\text{fo} = ({\rm B^\text{fo}} {\rm h^\text{fo}} - 2 t^\text{fo} ({\rm h^\text{fo}}-t^\text{fo}))
\end{equation}
\begin{equation}
\begin{split}
v^\text{fo}_\text{tot} = &\, (2{\rm n^{\text{rm}}}(l_\text{y}^\text{wa}-2t^\text{fo}) + ({\rm n^{\text{rm}}} + 1)l_\text{x}^\text{wa})\cdot\\
&\,({\rm B^\text{fo}} {\rm h^\text{fo}}-2t^\text{fo}({\rm h^\text{fo}}-t^\text{fo}))
\end{split}
\end{equation}
\begin{equation}
F^\text{fo}_\text{D,x} = P^\text{fo}_\text{D} l_\text{x}^\text{wa} + F_\text{D,x}^\text{wa}
\end{equation}
\begin{equation}
e = 0.5 \cdot ({\rm B^\text{fo}}-2 t^\text{fo} - t^\text{wa})
\end{equation}
\begin{equation}
\sum_{w \in \mathcal{W}} {\rm \underline{t}^\text{wa}_{w}} x_w \leq  t^\text{wa} \leq {\rm \bar{t}^\text{wa}}
\end{equation}
\begin{equation}
P_\text{D}^\text{wa} = {\rm g} t^\text{wa} h^\text{wa}\sum_{w \in \mathcal{W}}\rho_w x_w
\end{equation}
\begin{equation}
F_1 = {\rm g}{\rm C_\text{d}}v^\text{wa} \sum_{w \in \mathcal{W}} \rho_w x_w
\end{equation}
\begin{equation}
q^\text{ro} = v^\text{slc} \cdot n^\text{slc} \cdot \sum_{r \in \mathcal{R}}(x_r\rho_r  \quad + \quad  {\rm R^{\text{co}}}\sum_{c \in \mathcal{C}}(x_c \rho_c ))
\end{equation}
\begin{equation}
\sum_{f \in \mathcal{F}} {\rm \underline{t}^\text{fo}_{f}} x_f \leq 0.5{\rm B^\text{fo}}
\end{equation}
\begin{equation}
P^\text{fo}_\text{D} = {\rm g}\cdot A^\text{fo} \cdot \sum_{f \in \mathcal{F}} x_f\rho_f
\end{equation}
\begin{equation}
\sum_{w \in \mathcal{W}} x_w =1 ,  x_w \in \{0,1\}
\end{equation}
\begin{equation}
\sum_{f \in \mathcal{F}} x_f =1 ,  x_f \in \{0,1\}
\end{equation}
\begin{equation}
\sum_{r \in \mathcal{R}} x_r =1 ,  x_r \in \{0,1\}
\end{equation}
\begin{equation}
\sum_{c \in \mathcal{C}} x_c =1 ,  x_c \in \{0,1\}
\end{equation}
\begin{equation}
{\rm \underline{h}^\text{wa}} \leq h^\text{wa} \leq {\rm \bar{h}^\text{wa}}
\end{equation}
\begin{equation}
l_\text{x}^\text{fl},l_\text{y}^\text{fl} \in \{\underline{l}^\text{fl},\bar{l}^\text{fl}\} 
\end{equation}
\begin{equation}
n^\text{re} \geq {\rm \underline{n}^\text{re}}  , n^\text{re} \in \mathbb{Z}_+
\end{equation}
\begin{equation}
w^\text{do} \geq {\rm \underline{w}^\text{do}}
\end{equation}
\begin{equation}
l^\text{wi} \geq {\rm \underline{l}^\text{wi}}
\end{equation}
\begin{equation}
n^\text{slc} \leq {\rm \bar{n}^\text{slc}}, n^\text{slc} \in \mathbb{Z}_+
\end{equation}
\begin{equation}
x^e \in \{ 0,1\}
\end{equation}
\begin{equation}\label{eq:boundn:smooth}
x^{\text{wa}} \in \{ 0,1\}
\end{equation}
\end{small}

\section{Computing times}
\label{app_b}
For each of the Pareto fronts in this paper, Table~\ref{tab:soltime} lists the number of MINLP instances that were solved by the $\epsilon$-constraint method, the total computation time (in seconds), and the average time required for one MINLP instance.

\begin{table}[h]
    \caption{Solution time (in seconds) for Pareto-front problems}
    \label{tab:soltime}
\begin{tabular}{l r r l}
\toprule
  Problem &	Number of instances  & Solution time & AST$^*$ \\
  \midrule
All material available (Fig. \ref{fig:PF-3room-allMat})&	149	& 267.935	& 1.798\\
No soil blocks available (Fig. \ref{fig:MaterialSelection-3-room-exceptsoil})&	182 &	323.594	&1.778\\
Fixed wall material(brick G1)(Fig. \ref{fig:pareto-wallmat-allinall})&	19	&25.153	&1.324\\
Fixed wall material(brick G2)(Fig. \ref{fig:pareto-wallmat-allinall})&	169	&305.047	&1.805\\
Fixed wall material(cement G1)(Fig. \ref{fig:pareto-wallmat-allinall})&	6	&7.676	&1.279\\
Fixed wall material(cement G2)(Fig. \ref{fig:pareto-wallmat-allinall})&	14	&32.298	&2.307\\
Fixed wall material(soil G1)(Fig. \ref{fig:pareto-wallmat-allinall})&	12	&18.500	&1.542\\
Fixed wall material(soil G2)(Fig. \ref{fig:pareto-wallmat-allinall})&	12 &	29.460	&2.455\\
Varying floor area (Fig. \ref{fig:PF-area})&	134	& 238.484	&1.780\\
${\rm B}^\text{fo} =0.81$ (Fig.~\ref{fig:PF-0.81}.a)	&272 &	483.100	&1.776\\ 
${\rm B}^\text{fo} =0.81$, no stone G2 (Fig.~\ref{fig:PF-0.81}.b) &	279	& 498.072 &	1.785\\
\bottomrule
{$^*$ AST: average solution time per instance}
\end{tabular}
\end{table}
\end{appendices}
\bibliography{paper_finalazed}
\end{document}